\documentclass[12pt,reqno]{article}
\usepackage{amssymb}
\usepackage{amsmath}
\usepackage{mathrsfs}
\usepackage{amsthm}
\usepackage{amscd}
\usepackage{amsbsy}
\usepackage{amsfonts}
\usepackage{appendix}
\usepackage{float}
\usepackage{xcolor}
\usepackage{epstopdf}
\usepackage{graphicx}
\usepackage{psfrag}
\usepackage{subfigure}
\usepackage{overpic}
\usepackage{longtable}
\usepackage{multirow}
\usepackage{verbatim}
\usepackage{listings}
\usepackage{epsfig}
\usepackage{enumerate}
\usepackage{booktabs}
\usepackage{bm}
\usepackage{stackengine}
\newtheorem{thm}{Theorem}[section]

\newtheorem{cor}[thm]{Corollary}
\newtheorem{lem}[thm]{Lemma}
\newtheorem{prop}[thm]{Proposition}

\newtheorem{rem}[thm]{Remark}

\newcommand\xrowht[2][0]{\addstackgap[.5\dimexpr#2\relax]{\vphantom{#1}}}
\renewcommand{\theequation}{\thesection.\arabic{equation}}

\def\N{{\mathbb N}}

\def\R{{\mathbb R}}

\def\bb{\begin}
\def\bc{\begin{center}}       \def\ec{\end{center}}
\def\be{\begin{equation}}     \def\ee{\end{equation}}
\def\ba{\begin{array}}        \def\ea{\end{array}}
\def\bea{\begin{eqnarray}}    \def\eea{\end{eqnarray}}
\def\beaa{\begin{eqnarray*}}  \def\eeaa{\end{eqnarray*}}
\def\hh{\!\!\!\!}             
\def\EQ{\hh & = & \hh}        
      \def\GE{\hh & \ge & \hh}

\def\e{\varepsilon}           
               
\def\la{\lambda}

\def\oo{\infty}               \def\cd{\cdot}
\def\d{\cdot}
\def\dd{\cdots}               \def\pa{\partial}
\def\q{\quad}                 \def\qq{\qquad}
\def\f{\frac}                 
\def\z{\left}                 \def\y{\right}

\def\Lp{{\mathcal L}^p}

\def\nlp{\|\cdot\|_p}

\def\nmp#1{\|#1\|_p}

\def\nmt#1{\|#1\|_2}

\def\bpr{B_p[r]}

\def\pps{\phi_{p^*}}

\def\dx{\,{\rm d}x}

\def\rd{\,{\rm d}}

\def\lb{\label}
\def\x#1{{\rm (\ref{#1})}}

\def\andq{\qquad \mbox{and} \qquad}
\def\ifl{\iffalse}
\def\Proof{\noindent{\bf Proof} \quad}
\def\qed{\hfill $\Box$ \smallskip}

\def\qqf{\qquad \forall}

\begin{document}

\title{\bf On the Meromorphic Integrability of
the Critical Systems for Optimal Sums of Eigenvalues}

\author{Yuzhou Tian$^{a,}$\footnote{Correspondence author.},\qq Meirong Zhang$^{b}$}

\date{\today}%

\maketitle

\begin{center}
$^a$ Department of Mathematics, Jinan University, Guangzhou 510632, China	\\
$^b$ Department of Mathematical Sciences, Tsinghua University, Beijing 100084, China\\
E-mail: {\tt
tianyuzhou2016@163.com (Y. Tian)}\\
E-mail: {\tt zhangmr@tsinghua.edu.cn (M. Zhang)}
\end{center}

\begin{abstract}
The popularity of estimation to bounds for sums of eigenvalues started from P. Li and S. T. Yau for the study of the  P\'{o}lya conjecture. This subject is extended to different types of differential operators. This paper explores for  the sums of the first $m$  eigenvalues of Sturm-Liouville operators from two aspects.

Firstly, by the complete continuity of eigenvalues, we propose a family of critical systems consisting of nonlinear ordinary differential equations, indexed by the exponent $p\in(1,\infty)$ of the Lebesgue spaces concerned. There have profound relations between the solvability of these systems and the optimal lower or upper bounds for the sums of the first $m$  eigenvalues of Sturm-Liouville operators, which provides a novel idea to study the optimal bounds.

Secondly, we investigate the integrability or solvability of the critical systems. With suitable selection of exponents $p$, the critical systems are equivalent to the polynomial Hamiltonian systems of $m$ degrees of freedom.  Using the differential Galois theory, we perform a complete classification  for meromorphic integrability of these polynomial critical systems.  As a by-product of this classification, it gives a positive answer to the conjecture raised by Tian, Wei and Zhang [J. Math. Phys. 64, 092701 (2023)] on  the critical systems for optimal eigenvalue gaps. The numerical simulations  of the Poincar\'{e} cross sections show that the critical systems for sums of eigenvalues can appear complex dynamical phenomena, such as periodic trajectories, quasi-periodic trajectories and chaos.
\end{abstract}

{\small {\bf Mathematics Subject Classification (2020)}:
Primary 34L15; 
Secondary 37J30, 
70H07. 

{\bf Keywords:} Sums of eigenvalues, Sturm-Liouville operator, Critical system, Meromorphic integrable, Differential Galois theory.}


\section{Introduction} \lb{mr}

In 1807, the pioneering work of solving the heat equation by Fourier  planted the seed for the spectral theory of differential operators. Inspired by Fourier's work, Sturm and
Liouville in 1837 systematically treated the spectra of second-order linear ordinary differential operators, commonly referred to as Sturm-Liouville operators.
 Afterwards, their work gradually evolved a whole new branch of mathematics, namely Sturm-Liouville theory. In the 20th century, Weyl's famous work \cite{MR1511560} together with the birth of quantum mechanics  revolutionize this theory.
 Henceforth, the modern Sturm-Liouville theory not only provides a perfect medium for understanding the quantum mechanics, but also greatly promotes the development of other areas of mathematics, such as harmonic analysis, differential geometry and operator algebras.  Nowadays, Sturm-Liouville theory is still an active area of research in modern mathematical physics.

 Let $\Omega=\left[0,1\right]$
  be the unit interval. Fixed an exponent $p\in\left(1,\infty\right)$, the $L^p$ Lebesgue space on $\Omega$ is denoted by
 $$\Lp:=L^p\left(\Omega,\R\right).$$
 For an integrable potential $q\in\Lp$, we consider the following Dirichlet eigenvalue problem for the \emph{Sturm-Liouville operator or
one-dimensional Schr\"{o}dinger operator}
\begin{equation}\label{line}
	\begin{split}
\mathscr{D}_q \psi &:=-\psi'' + q \psi=\la \psi, \qq x\in \Omega,\\
&\psi\big|_{\partial\Omega}=0.
	\end{split}
\end{equation}
A number $\la$ is an \emph{eigenvalue} of the system \eqref{line} if it has a nontrivial solution $\psi\left(x\right)$, called an \emph{eigenfunction} associated $\la$. It is
well-known that
the eigenvalues of problem  \eqref{line}
can be written in the form of an increasing sequence
 \[
 \la_1(q) < \la_2(q) < \dd < \la_m(q) < \dd, \qq \la_m(q) \to +\oo \mbox{ as } m\to \oo.
 \]
 Here we have regarded  eigenvalues as nonlinear functionals of potentials $q\in \Lp$. The sum of the first $m$ eigenvalues is defined as
 $$\mathscr{E}_m\left(q\right):=\sum_{i=1}^m\lambda_i\left(q\right).$$

 In quantum theory, the eigenvalues have definite physical significance, which  correspond to the energy levels of a particle within a potential $q$.  Thereby, $\mathscr{E}_m\left(q\right)$ is the \emph{total energy} of $m$ particles. 
 Especially,  the absorption energy for particle from the ground state to the first excited state is described by the \emph{fundamental eigenvalue gap} $\la_2\left(q\right)-\la_1\left(q\right)$. Because of the above physical interpretations, estimation to the lower and upper bounds  for eigenvalue problems, including gap, ratio and sum, are central to a large part of Sturm-Liouville theory.

For the lower bounds of the fundamental eigenvalue gaps, many fascinating results about different types of operators and boundary conditions have been contributed by a lot of  mathematicians up to the present days, see for example \cite{MR2784332,MR3478937,MR2299195,MR2511892,MR4150221,MR2881964,MR1948113,MR1614731,MR1081670} and references therein. The estimate for the upper bounds can be found in \cite{MR4447102,MR2881964,MR1948113,MR1614731,MR829055,MR1081670}.
The problems of estimations for the eigenvalue ratios also have been extremely studied in
 \cite{MR1218744,MR257592,MR4339006}.

In applied sciences and mathematics, it is important to understand the sums of eigenvalues.  For example, related with the quantum mechanics,  elasticity theory, geometry and PDEs, it
is natural to study the sum of the first $m$ eigenvalues \cite{MR639355,MR4379307}.  Perhaps the most important motivation is to originate from the famous P\'{o}lya conjecture \cite{MR129219} about the lower bound on $i$-th eigenvalue for the Dirichlet Laplacian, which still
remains open.  In 1983, Li and Yau  \cite{MR701919} gave a partial answer to this conjecture
and  improved Lieb's result \cite{MR573436}. In order to be close to  P\'{o}lya conjecture, their technique is to estimate the lower bound for the sum of the first $m$ eigenvalues, commonly known as the \emph{Berezin-Li-Yau bound or  inequality}.  The estimates for sum of eigenvalues of different types of operators have gained wide investigation  from mathematicians since Li and Yau. We refer the readers to \cite{MR2782621,MR1165859,MR1383015,MR4379307,MR2846268,MR3412394,MR3488539}, etc. But up to now,  there is not a general method to obtain the optimal
lower or upper bound
for sum of  eigenvalues of Sturm-Liouville operator \eqref{line}.

The purpose of this work is to investigate the optimization problems on sum of eigenvalues for \eqref{line}.  Let
$$B_{p,r}=\left\{q\in \Lp:\;\parallel q\parallel_p\leq r\right\}$$
be the (infinitely dimensional) ball of radius $r$, centered at the origin, in the space $\left(\Lp,\parallel\cdot\parallel_p\right)$.
We consider  the following optimization problems
\begin{align}\label{eq66}
&\mathscr{E}_m^-:=\min\limits_{q\in B_{p,r}}\mathscr{E}_m\left(q\right)\;\text{and}\; \mathscr{E}_m^+:=\max\limits_{q\in B_{p,r}}\mathscr{E}_m\left(q\right).
\end{align}
Their solutions will provide the following estimations on sum of the eigenvalues
\begin{align}\label{eq67}
&\mathscr{E}_m^-\leq\mathscr{E}_m\left(q\right)\leq\mathscr{E}_m^+, \qqf q\in \Lp.
\end{align}
The lower bound $\mathscr{E}_m^-$ is also called \emph{Berezin-Li-Yau type lower bound}. Remarkably,  $\mathscr{E}_m^-$ and $\mathscr{E}_m^+$ are the optimal lower and upper bounds of  $\mathscr{E}_m\left(q\right)$ in a certain sense, respectively.

Our first result provides a completely different  approach to attain possible solutions to problems \eqref{eq66}. As a consequence of the complete continuity of eigenvalues in potentials \cite{MZ10, YZ11, Zh08}, one shows that the optimization problems \eqref{eq66}  can be attained by some optimizing potentials $q^{\pm}\in B_{p,r}$. See Theorem \ref{attain}.

In order to determine $q^\pm$ and $\mathscr{E}_m^\pm$, we establish the next result.
\begin{thm}\label{main}
	Let the exponent $p\in(1,\oo)$, $r\in(0,\oo)$ and $m\in\mathbb{N}$ be given with $m\geq2$. Denote by $p^*:=p/(p-1) \in (1,\oo)$  the conjugate exponent of $p$. For problems \eqref{eq66}, indicated by $\e=-$ and $+$ respectively,  there are $m$-dimensional parameters $\left(\mu_1,\ldots,\mu_m\right)=\left(\mu_1^\e,\ldots,\mu_m^\e\right)$ and non-trivial solutions $(u_1(x),\ldots,u_m(x))=(u_1^\e(x),\ldots,u_m^\e(x))$ to the following system
	\begin{align}\label{ceq37}
		&-u''_i+\e\left(\sum_{j=1}^mu_j^2\right)^{p^*-1}u_i=\mu_iu_i,\quad i=1,\ldots,m,
	\end{align}
	such that
	
	\emph{(i)} the solutions  $u_i(x)$ satisfy the Dirichlet boundary condition for $i=1,\ldots,m$.
	
	\emph{(ii)} the solutions $u_i(x)$ satisfy
	\be \lb{u12}
	\int_\Omega \left(\sum_{i=1}^mu_i^2(x)\right)^{p^*} \dx = r^p,
	\ee
	and the optimizing potentials $q^\e(x)$ are determined by
	\be \lb{qx2}
	q^\e(x):= \e\z(\sum_{i=1}^m\left(u_i^\e(x)\right)^2\y)^{p^*-1}, \qq x\in \Omega.
	\ee
	
	\emph{(iii)} the minimal and maximal of the sum of the first $m$ eigenvalues  are given by
	\be \lb{LM11}
\mathscr{E}_m^-=\sum_{i=1}^m \mu_i^-\;\text{and}\; \mathscr{E}_m^+=\sum_{i=1}^m \mu_i^+,
	\ee
respectively.
\end{thm}

 System \eqref{ceq37} is called in this paper {\it the critical system}, which is deduced by the direct application of the Lagrangian multiplier method to problems \eqref{eq66}, as done in \cite{WMZ09, Zh09,YZ12}. Compared with the deductions of the critical systems to the inverse spectral problems for elliptic operators or Sturm-Liouville operators by Ilyasov and Valeev  \cite{MR3912726,MR4197914,MR4285920}, our derivation approach, employed  the complete continuity of eigenvalues $\lambda_m\left(q\right)$ in $q\in \Lp$ \cite{MZ10, Zh08}, is very simpler.

Let $v_i=u'_i$ for $i=1,\ldots,m$.  Then critical system  \eqref{ceq37} is equivalent to a Hamiltonian system of $m$ degrees of freedom
\begin{align}\label{ceq38}
	&u'_i=v_i=\frac{\partial H}{\partial v_i},\quad v'_i=-\mu_iu_i+\varepsilon\left(\sum_{j=1}^mu_j^2\right)^{p^*-1}u_i=-\frac{\partial H}{\partial u_i},\quad  i=1,\ldots,m,
\end{align}
with the Hamiltonian
\begin{align}\label{ceq39}
	&H=\frac{1}{2}\sum_{i=1}^m\left(v_i^2+\mu_iu_i^2\right)-\frac{\varepsilon}{2p^*}\left(\sum_{j=1}^mu_j^2\right)^{p^*}.
\end{align}

Let  $\mathbf{u}=\left(u_1,\ldots,u_m\right)$ and  $\mathbf{v}=\left(v_1,\ldots,v_m\right)$. The non-constant function $I=I\left(\mathbf{u},\mathbf{v }\right)$ is said to be a \emph {first integral} of the Hamiltonian system \eqref{ceq38} if $H$ and $I$ are \emph {in involution}, i.e.
the Poisson bracket
\begin{align}\label{PB}
	&\left\{H,I\right\}=\sum_{i=1}^m\left(\frac{\partial H}{\partial v_i}\frac{\partial I}{\partial u_i}-\frac{\partial H}{\partial u_i}\frac{\partial I}{\partial v_i}\right)=0.
\end{align}
The Hamiltonian function  $H$ itself is always a first integral due to the  antisymmetry of Poisson bracket. Denote the gradient of function $I$ by $\nabla I$.
The functions $I_i$ for $i=1,\ldots,l$ are  \emph{functionally independent} on $U$ if
$$\text{rank}\left(\nabla I_1,\ldots,\nabla I_l\right)=l$$
with the possible exception
of sets of Lebesgue measure zero. The Hamiltonian system \eqref{ceq38} is called \emph{completely integrable}, or simply \emph{integrable}  in Liouville's sense if there exist $m$  functionally independent first integrals $I_1\equiv H,I_2,\ldots,I_m$ ($H$ is the Hamiltonian). In addition, Hamiltonian system \eqref{ceq38} is \emph{meromorphic completely integrable},  or simply \emph{meromorphic integrable} if its $m$ functionally independent first integrals $I_1\equiv H,I_2,\ldots,I_m$ are meromorphic.

Theorem \ref{main} allows us to determine a solution to  the optimization problems \eqref{eq66} by solving a boundary value problem for critical system  \eqref{ceq37}. In other words, the solvability of Hamiltonian system \eqref{ceq38} means  the solvability of problems \eqref{eq66}.  The classical Arnold-Liouville's theorem  exhibits that if Hamiltonian system \eqref{ceq38} is completely integrable, then it can be solved by quadrature, see \cite{MR2004534}. Naturally, we will focus on the next problem.

\vspace{2ex}
\noindent{\bf Problem 1.}  \emph{Whether or not Hamiltonian system \eqref{ceq38} is completely integrable.}
\vspace{2ex}

\noindent The answer is too difficult because of the two reasons: system  \eqref{ceq38} is not a  polynomial system  for some $p^*\in\left(1,+\oo\right)$; there are no universal techniques  to decide the integrability of Hamiltonian systems.	 On the other hand, of particular interest is the limiting case of problems \eqref{eq66}, that is, $p=1$. For this limiting case,  like in \cite{WMZ09, YZ12, Zh09}, it is significant to investigate the limiting system of \x{ceq37} as $p\downarrow 1$, i.e., as $p^* \uparrow \oo$. In such a limiting process, let us pay a special attention to  exponents $p$ so that
\begin{align}\label{eq70}
p= \frac{k}{k-1}\;\text{and}\;
p^*=k, \;\text{and}\;k=2,3, \dd
\end{align}
For these exponents, system  \eqref{ceq38} are the following  polynomial Hamiltonian systems
\begin{align}\label{eq38}
	&u'_i=v_i,\quad v'_i=-\mu_iu_i+\varepsilon\left(\sum_{j=1}^mu_j^2\right)^{k-1}u_i,\quad  i=1,\ldots,m,
\end{align}
where
\begin{align}\label{eq39}
	&H=\frac{1}{2}\sum_{i=1}^m\left(v_i^2+\mu_iu_i^2\right)-\frac{\varepsilon}{2k}\left(\sum_{j=1}^mu_j^2\right)^k.
\end{align}

At present, most studies have been dedicated to the integrability for Hamiltonian system of  two degrees of freedom, see for instance \cite{MR3743739,MR943701,MR4075569,MR2831794} and references therein. Except the natural Hamiltonian system with homogeneous potential, there are few literature relating to the integrability of other types of Hamiltonian systems of arbitrary degrees of freedom, see \cite{MR2123446, MR1764944,MR3804721}. With the help of the differential Galois theory, we give a complete classification of meromorphic integrability for Hamiltonian system \eqref{eq38} as follows.

\begin{thm}\label{th1}
	Let  $\mathscr{U}=\left(\mu_1,\ldots,\mu_m\right)$ and $k\in \mathbb{N}^+$ with $k\geq 2$. The Hamiltonian system \eqref{eq38} is meromorphic completely integrable if and only if $k$ and $\mathscr{U}$ belong to one of the following two families:
	\begin{table}[H]
		\centering
		\caption{Meromorphic completely integrable cases}
		\begin{tabular}{|c|c|c|c|}
			\hline
			Case&$k$&$\mathscr{U}$&Additional meromorphic first integrals\xrowht{20pt}\\
			\hline
			1&$k=2$&$\mathscr{U}\in \mathbb{R}^m$& See Proposition \ref{pr5}.\xrowht{28pt}\\
			\hline
			2&$k\geq3$&$\mu_1=\mu_2=\cdots=\mu_m$&$I_i=u_1v_{i+1}-u_{i+1}v_1,\quad i=1,\ldots,m-1$.\xrowht{28pt}\\
			\hline
		\end{tabular}%
		\label{Ta2}%
	\end{table}%
\end{thm}



Especially when $\varepsilon=-$, $m=2$ and $k=2n$, the next corollary can be attained by  the linear canonical transformation $\left(u_1,u_2,v_1,v_2\right)\mapsto \left(-u_1,-\text{i}u_2,-v_1,\text{i}v_2\right)$.
\begin{cor}\label{co1}
Consider the following	Hamiltonian system
\begin{equation}\label{e12}
	\begin{split}
	u'_1=v_1,\quad	u'_2=-v_2,\quad v'_1=-\mu_1u_1-u_1\left(u_1^2-u_2^2\right)^{2n-1},\quad v'_2=\mu_2u_2+u_2\left(u_1^2-u_2^2\right)^{2n-1}
	\end{split}
\end{equation}
	with Hamiltonian
	\begin{align}\label{e15}
		&H= \f{1}{2}\z(v_1^2-v_2^2 + \mu_1 u_1^2 - \mu_2 u_2^2\y) +\f{1}{4n}\z(u_1^2-u_2^2\y)^{2n}.
	\end{align}
For $\mu_1\neq \mu_2$, $n\in \mathbb{N}^+$ and $n\geq2$, system \eqref{e12}
is meromorphic non-integrable.
\end{cor}

In \cite{TWZ}, the authors studied the critical system for optimal eigenvalue gaps and posed the next conjecture.

\vspace{2ex}
\noindent{\bf Conjecture}. For generic $\mu_1\neq \mu_2$  and $n\geq2$,
system \eqref{e12} is not polynomial integrable.
\vspace{2ex}

Obviously, Corollary \ref{co1} not only gives a positive answer to  the above conjecture, but also  extends it  to meromorphic non-integrable.

The framework of the paper is as follows. We briefly recall some preliminary concepts and results of differential Galois approach in section \ref{se2}.  After  gathering the complete continuity results on eigenvalues, we deduce the critical system \eqref{ceq37} in section \ref{se3}. To prove Theorem \ref{th1}, we will divide into two sections. In section  \ref{se4}, we show that system \eqref{eq38} is complete integrability if the parameters $k$ and $\mathscr{U}$ belong to Table \ref{Ta2}.  In section \ref{se5}, by Morales-Ramis theory, we prove  that system \eqref{eq38} is meromorphic non-integrability when the parameters $k$ and $\mathscr{U}$ are outside Table \ref{Ta2}. Section \ref{se6} presents exemplary Poincar\'{e} cross sections of the critical system \eqref{eq38}, which exhibits that system \eqref{eq38} has abundant dynamical behaviors.

\section{Preliminaries}\label{se2}
In this section, we introduce some necessary concepts and preliminary results, containing Morales-Ramis theory, Hypergeometric equation and Kovacic's results.
\subsection{Morales-Ramis theory}
The Morales-Ramis theory  \cite{MR1713573} is a powerful tool to determine the non-integrability of complex Hamiltonian systems. Roughly speaking, this theory establishes a relation between the meromorphic  integrability and the differential
Galois group of the variational equations or the normal variational equations. Next we briefly describe the Morales-Ramis theory. For some precise notions of differential Galois theory, see \cite{MR1960772}. 

Consider a complex symplectic manifold $M\subset\mathbb{C}^{2m}$ of dimension $2m$  with the standard symplectic form $\bm{\tilde{\omega}}=\sum_{j=1}^mdu_j\wedge dv_j$.  Let $H: M\rightarrow \mathbb{C}$ be a holomorphic  Hamiltonian. The Hamiltonian system with $m$ degrees of freedom is given by
\begin{align}\label{eq59}
&\dfrac{d \mathbf{x}}{dt}=X_H\left(\mathbf{x}\right)=\left(\dfrac{\partial H}{\partial \mathbf{v}}, -\dfrac{\partial H}{\partial \mathbf{u}}\right),\quad t\in\mathbb{C}, \quad \mathbf{x}=\left(\mathbf{u},\mathbf{v}\right)\in M,
\end{align}
where $\mathbf{u}=\left(u_1,\ldots,u_m\right)$ and $\mathbf{v}=\left(v_1,\ldots,v_m\right)$ are the canonical coordinates.

Let $\Gamma$  be  a non-equilibrium solution of system  \eqref{eq59}. Assume that $\Gamma$ can be parameterized by time $t$, that is,
\begin{align*}
\bm{\varphi}:\mathbb{C}&\rightarrow M\subset\mathbb{C}^{2m}\\
t&\mapsto\left(\mathbf{u}\left(t\right),\mathbf{v}\left(t\right)\right).
\end{align*}
Then the \emph{variational equation} (VE) along $\Gamma$ is the linear differential system
\begin{align}\label{eq60}
&\dfrac{d \mathbf{y}}{dt}=\dfrac{\partial X_H\left(\bm{\varphi}\left(t\right)\right)}{\partial \mathbf{x}}\mathbf{y}, \quad \mathbf{y}\in T_\Gamma M,
\end{align}
where $T_\Gamma M$ is the tangent bundle $TM$ restricted on $\Gamma$.

Let $N:=T_\Gamma M/T\Gamma$ be the normal bundle
 of $\Gamma$ \cite{MR1411677}, and $\pi: T_\Gamma M \rightarrow N$ be the nature projective homomorphism. The \emph{normal variational equation} (NVE) along $\Gamma$ has the form
\begin{align}\label{eq61}
&\dfrac{d \mathbf{z}}{dt}=\pi_*\left(T\left(\mathfrak{u}\right)\left(\pi^{-1}\mathbf{z}\right)\right),\quad \mathbf{z}\in N,
\end{align}
where $\mathfrak{u}=X_H\left(\mathbf{x}\right)$ with $\mathbf{x}\in M$, and $T\left(\mathfrak{u}\right) $ is the tangential
variation of $\mathfrak{u}$ along $\Gamma$, that is, $T\left(\mathfrak{u}\right)=\partial X_H/\partial \mathbf{x}$.  Note that the above NVE is a $2\left(m-1\right)$-dimensional linear differential system. We can employ a generalization of D'Alambert's method to get the NVE \eqref{eq61}, see \cite{MR1713573}.
Briefly speaking, we use the fact that $X_H\left(\bm{\varphi}\left(t\right)\right)$ is a solution of the VE \eqref{eq60} to reduce its dimension  by one. In effect, we typically restrict the equation \eqref{eq59} to the energy level $h=H\left(\bm{\varphi}\left(t\right)\right)$. Then the dimension of  the corresponding VE \eqref{eq60} also can be reduced.

Morales and Ramis \cite{MR1713573} proved the following classical theorem, which give a necessary condition for the integrability of Hamiltonian system \eqref{eq59} in the Liouville sense.
\begin{thm}\label{th5}{\rm (Morales-Ramis theorem, see \cite{MR1713573})} If Hamiltonian system \eqref{eq59} is meromorphically integrable in the Liouville
sense in a neighbourhood of a particular solution $\Gamma$, then the identity component of the Galois group of the \emph{NVE} \eqref{eq61} is Abelian.
\end{thm}

The next theorem tells us that the identity component of the differential Galois group is invariant under the covering.
\begin{thm}\label{th4}{\rm (\cite{MR1713573})}
Let $\mathcal{M}$ be a connected Riemann surface and $\nabla$ be a meromorphic connection over $\mathcal{M}$. 	Assume that $f:\mathcal{M}'\longrightarrow\mathcal{M}$ is a finite ramified covering of $\mathcal{M}$ by a connected Riemann surface  $\mathcal{M}'$.  Let $\nabla'=f^*\nabla$, i.e. the pull back of $\nabla$ by $f$. Then there exists a natural
injective homomorphism
$$\emph{Gal}\left(\nabla'\right)\rightarrow\emph{Gal}\left(\nabla\right)$$
of differential Galois groups which induces an isomorphism between their Lie
algebras.
\end{thm}

\subsection{Hypergeometric equation}
The hypergeometric equation is a  second order differential equation over the Riemann sphere $\mathbf{P}^1$ with three  regular singular points \cite{MR2682403,MR0178117}. Let us consider the following form of hypergeometric equation with three singular points at $z=0,1,\infty$
\begin{align}\label{hy}
	&\dfrac{d^2\zeta}{dz^2}+\left(\dfrac{1-\alpha-\tilde{\alpha}}{z}+\dfrac{1-\gamma-\tilde{\gamma}}{z-1}\right)\dfrac{d\zeta}{dz}+\left(\dfrac{\alpha\tilde{\alpha}}{z^2}+\dfrac{\gamma\tilde{\gamma}}{\left(z-1\right)^2}+\dfrac{\beta\tilde{\beta}-\alpha\tilde{\alpha}-\gamma\tilde{\gamma}}{z\left(z-1\right)}\right)\zeta=0,
\end{align}
where $\left(\alpha,\tilde{\alpha}\right)$, $\left(\gamma,\tilde{\gamma}\right)$ and $\left(\beta,\tilde{\beta}\right)$ are the exponents at the respective singular points, and meet the Fuchs relation
$\alpha+\tilde{\alpha}+\gamma+\tilde{\gamma}+\beta+\tilde{\beta}=1$.
The exponent differences can be defined as
$\varrho=\alpha-\tilde{\alpha}$, $\varsigma=\gamma-\tilde{\gamma}$ and $\tau=\beta-\tilde{\beta}$.

The following theorem goes back to Kimura \cite{MR277789}, whose gave
necessary and sufficient
conditions for solvability of the identity component of the differential Galois group of  \eqref{hy}.
\begin{table}[H]
	\centering
	\caption{Schwarz table with $l,s,\upsilon\in\mathbb{Z}.$}\vspace{2ex}
	\begin{tabular}{|c|c|c|c|c|}
		\hline
		$1$&$1/2+l$&$1/2+s$&Arbitrary complex number&\xrowht{6pt}\\
		\hline
		$2$&$1/2+l$&$1/3+s$&$1/3+\upsilon$&\xrowht{6pt}\\
		\hline
		$3$&$2/3+l$&$1/3+s$&$1/3+\upsilon$&$l+s+\upsilon$ even\xrowht{6pt}\\
		\hline
		$4$&$1/2+l$&$1/3+s$&$1/4+\upsilon$&\xrowht{6pt}\\
		\hline
		$5$&$2/3+l$&$1/4+s$&$1/4+\upsilon$&$l+s+\upsilon$ even\xrowht{6pt}\\
		\hline
		$6$&$1/2+l$&$1/3+s$&$1/5+\upsilon$&\xrowht{6pt}\\
		\hline
		$7$&$2/5+l$&$1/3+s$&$1/3+\upsilon$&$l+s+\upsilon$ even\xrowht{6pt}\\
		\hline
		$8$&$2/3+l$&$1/5+s$&$1/5+\upsilon$&$l+s+\upsilon$ even\xrowht{6pt}\\
		\hline
		$9$&$1/2+l$&$2/5+s$&$1/5+\upsilon$&$l+s+\upsilon$ even\xrowht{6pt}\\
		\hline
		$10$&$3/5+l$&$1/3+s$&$1/5+\upsilon$&$l+s+\upsilon$ even\xrowht{6pt}\\
		\hline
		$11$&$2/5+l$&$2/5+s$&$2/5+\upsilon$&$l+s+\upsilon$ even\xrowht{6pt}\\
		\hline
		$12$&$2/3+l$&$1/3+s$&$1/5+\upsilon$&$l+s+\upsilon$ even\xrowht{6pt}\\
		\hline	
		$13$&$4/5+l$&$1/5+s$&$1/5+\upsilon$&$l+s+\upsilon$ even\xrowht{6pt}\\
		\hline
		$14$&$1/2+l$&$2/5+s$&$1/3+\upsilon$&$l+s+\upsilon$ even\xrowht{6pt}\\
		\hline		
		$15$&$3/5+l$&$2/5+s$&$1/3+\upsilon$&$l+s+\upsilon$ even\xrowht{6pt}\\
		\hline							
	\end{tabular}%
	\label{Ta3}%
\end{table}%
\begin{thm}\label{A1}{\rm (\cite{MR277789})} The identity component of the Galois group of the hypergeometric equation \eqref{hy} is solvable if and only if  either
		\begin{itemize}
		\item [\emph{(i)}] at least one of the four numbers    $\varrho+\tau+\varsigma,-\varrho+\tau+\varsigma,\varrho-\tau+\varsigma,\varrho+\tau-\varsigma$ is an odd integer, or
		\item [\emph{(ii)}] the numbers $\varrho$ or $-\varrho$, $\varsigma$ or $-\varsigma$ and $\tau$ or  $-\tau$ belong (in an arbitrary order)
		to some of the following fifteen families, see Table \ref{Ta3}.
	\end{itemize}
\end{thm}

\subsection{Kovacic's results}
Let $\mathbb{C}\left(z\right)$ be the field of rational functions  in the variable $z$ with complex coefficients. Consider the second order linear differential equation
\begin{align}\label{eq28}
	&\chi''=r\left(z\right)\chi,\quad r\left(z\right) \in \mathbb{C}\left(z\right).
\end{align}
It is well known that  the differential Galois group $G$ of equation \eqref{eq28} is an algebraic subgroup of $\text{SL}\left(2,\mathbb{C}\right)$.  In 1986, Kovacic \cite{MR839134}
characterized all possible types of $G$ as follows.

\begin{thm}\label{th6}{\rm (\cite{MR839134})}
 The differential Galois group $G$ of equation \eqref{eq28} is conjugated to one of the following four types:
		\begin{itemize}
		\item [\emph{(i)}] $G$ is conjugated to a subgroup of a triangular group, and equation \eqref{eq28}  admits a solution of the form $\chi=\exp\left(\int \omega \right)$ with $\omega \in\mathbb{C}\left(z\right)$.
		\item [\emph{(ii)}] $G$  is conjugate to a subgroup of
		\begin{align*}
			&\mathcal{G}=\left\{
			\left(
			\begin{array}{cc}
				\mathfrak{a}&0\\
				0&\mathfrak{a}^{-1}\\
			\end{array}
			\right)
			\Bigg|\mathfrak{a}\in \mathbb{C}\setminus\left\{0\right\}\right\}
			\bigcup
			\left\{
			\left(
			\begin{array}{cc}
				0&\mathfrak{a}\\
				\mathfrak{a}^{-1}&0\\
			\end{array}
			\right)
			\Bigg|\mathfrak{a}\in \mathbb{C}\setminus\left\{0\right\}
			\right\},
		\end{align*}
	and equation \eqref{eq28}  admits a solution of the form $\chi=\exp\left(\int \omega \right)$, where $\omega$ is algebraic of degree $2$ over $\mathbb{C}\left(z\right)$.
		\item [\emph{(iii)}] $G$ is finite and all solutions of equation \eqref{eq28} are algebraic over $\mathbb{C}\left(z\right)$.
		\item [\emph{(iv)}] $G=\emph{SL}\left(2,\mathbb{C}\right)$ and equation \eqref{eq28} does not admit Liouvillian solution.
	\end{itemize}
\end{thm}

Let $r\left(z\right)=\mathfrak{p}\left(z\right)/\mathfrak{q}\left(z\right)$
with $\mathfrak{p}\left(z\right),\mathfrak{q}\left(z\right)\in\mathbb{C}\left[z\right]$ relatively prime. The \emph{pole} of $r\left(z\right)$ is a zero of $\mathfrak{q}\left(z\right)$ and \emph{the order of the pole} is the multiplicity of the zero of $\mathfrak{q}\left(z\right)$.  \emph{The order of $r\left(z\right)$ at $\infty$} is defined by $\deg \mathfrak{q}-\deg \mathfrak{p}$.
Kovacic \cite{MR839134} also
provided the necessary conditions for types (i), (ii), or (iii) in  Theorem \ref{th6} to occur.
\begin{prop}\label{pr1}{\rm (\cite{MR839134})}
For the first three types in Theorem \ref{th6}, the necessary conditions of occurrence are respectively as follows:
\begin{description}
\item[\emph{\bf Type (i)}] Each pole of $r\left(z\right)$ must have even order or else have order $1$. The order of $r\left(z\right)$ at $\infty$ must be
even or else be greater than $2$.
\item[\emph{\bf Type (ii)}] The rational function $r\left(z\right)$ must have at least one pole that either has odd order greater than $2$ or else has
order $2$.
\item[\emph{\bf Type (iii)}] The order of a pole of $r\left(z\right)$ cannot exceed $2$ and the order of $r\left(z\right)$ at $\infty$ must be at
least $2$. If the partial fraction decomposition of $r\left(z\right)$ is
$$r\left(z\right)=\sum_i\dfrac{\alpha_i}{\left(z-c_i\right)^2}+\sum_j\dfrac{\beta_j}{z-b_j},$$
then $\sqrt{1+4\alpha_i}\in\mathbb{Q}$ for each $i$, $\sum_j\beta_j=0$, and if $\Delta=\sum_i\alpha_i+\sum_j\beta_j$, then $\sqrt{1+4\Delta}\in\mathbb{Q}$.
\end{description}
\end{prop}

\begin{rem}\label{re1}
For a general second order linear differential equation
\begin{align*}
	&y''=a_1y'+a_2,\quad a_1,a_2\in \mathbb{C}\left(z\right),
\end{align*}
 it can be transformed into the form \eqref{eq28} with
 $$r\left(z\right)=\dfrac{a_1^2}{4}-\dfrac{a'_1}{2}+a_2$$
 via the change
 \begin{align}\label{eq63}
 	&y=\exp\left(\dfrac{1}{2}\int a_1 dz\right)\chi.
 \end{align}
\end{rem}

\section{Deduction of the critical systems}\label{se3}

To derive the critical system \eqref{ceq37},  we refer to some basic properties of eigenvalues, viewed as functionals of potentials.

\bb{lem}{\rm (\cite{PT87,YZ11})}\lb{evsd}
Given $m\in\N$, the $m$th eigenvalue $\la_m(q)$ is continuously Fr\'echet differentiable in $q\in
(\Lp,\nlp)$, $p\in[1,\oo]$. Moreover, the Fr\'echet derivative $\pa_q \la_m(q)$, considered as an element of the conjugate space
$(\Lp)^*$, is
\be \lb{lamd}
\pa_q \la_m(q) =  (E_m(\cd;q))^2,
\ee
where $E_m(x)=E_m(x;q)$ is an eigenfunction associated with $\la_m(q)$ satisfying the following normalization condition
\be \lb{lamd1}
\nmt{E_m}= \z(\int_\Omega E_m^2(x)\dx\y)^{1/2}=1, \andq E'_m(0)>0.
\ee
\end{lem}

Let $q_l,\ q\in\Lp$ with $p\in[1,\oo]$. We say that $q_l$ is \emph{weakly
	convergent} to $q$ in $\Lp$ with respect to the weak topology $w_p$ if $$\lim\limits_{l\rightarrow\oo}\int_\Omega q_l\left(x\right)\xi\left(x\right)dx=\int_\Omega q\left(x\right)\xi\left(x\right)dx,\qquad\forall \xi\in\mathcal{L}^{p^*}.$$
Such a convergence is also denoted by $q_l \rightharpoonup q$ in $\Lp$.

The following lemma is the complete continuity of eigenvalues in weak topologies, see \cite{MZ10,Zh08,YZ11} for more details.

\bb{lem}{\rm (\cite{MZ10,Zh08,YZ11})}\lb{evsc}
Given $m\in\N$, the $m$th eigenvalue $\la_m(q)$ is completely continuous in $q\in(\Lp,w_p)$, $p\in[1,\oo]$. Here $w_p$ indicates the topology of weak convergence. More precisely, whenever $q_l \rightharpoonup q$ in $\Lp$, one has
\be \lb{lamc}
\lim_{l\to\oo} \la_m(q_l) = \la_m(q).
\ee
\end{lem}

\begin{thm}\label{attain}
Let $p\in (1,\oo)$, $r\in(0,\oo)$ and $m \in\mathbb{N}$ be given with $m\geq2$. Then there exist potentials $q^\e= q_{m,p,r}^\e\in \Lp$, $\e=+, \, -$ such that
\be \lb{ppm}
\nmp{q^-} =\nmp{q^+} =r,
\ee
and
\be \lb{ppm1}
\mathscr{E}_m^-=\mathscr{E}_m\left(q^-\right),\qq\mathscr{E}_m^+= \mathscr{E}_m^+\left(q^+\right).
\ee	
\end{thm}

\Proof
Let
$$B_{p,r}=\left\{q\in \Lp:\;\parallel q\parallel_p\leq r\right\}$$
be the ball of the space $(\Lp,\nlp)$. As we know, the ball   $B_{p,r}$ is a compact set of $(\Lp, w_p)$ with $p\in (1,\oo)$.  Since the sum of the first $m$ eigenvalues       $\mathscr{E}_m\left(q\right)$ is a finite sum, by the Lemma \ref{evsc}, then  $\mathscr{E}_m\left(q\right)$  is completely continuous in $q\in \Lp$. So, there exist $q^\pm(x)= q_{m,p,r}^\pm(x)\in \Lp$ such that
\(
\nmp{q^\pm} \le r
\)
and \x{ppm1} is confirmed.

Obviously, the Fr\'echet derivatives
$$\partial_q\left(\mathscr{E}_m\left(q\right)\right)\big|_{q=q^\pm}=\sum_{i=1}^m\left(E_i\left(\cd;q^\pm\right)\right)^2$$
is a non-zero function. From the Lagrange theorem,  it follows that the optimizing potentials $q^\pm$ cannot be such that $\nmp{q^\pm}<r$. This implies that $\nmp{q^\pm}=r$, that is, equation \eqref{ppm}.
\qed

Theorem \ref{attain}  tells us that problems \eqref{eq66} are constrained optimization problems, that is,
\begin{equation}\label{eq78}
\min \left(\max\right) \left(\mathscr{E}_m\left(q\right)\right)\;\text{subject to}\;\nmp{q} =r.
\end{equation}
Note that the Fr\'echet derivatives of  $\mathscr{E}_m\left(q\right)$ and the $L^p$ norm are
$$\partial_q\left(\mathscr{E}_m\left(q\right)\right)=\sum_{i=1}^m\left(E_i\left(x;q^\pm\right)\right)^2$$
and
 \[
\pa_q \nmp q = \nmp q^{1-p} |q(x)|^{p-2} q(x), \qq q\ne 0,
\]
respectively.

One can perform the Lagrange multiplier method to problems \eqref{eq78} to obtain that $q=q_{m,p,r}^\pm$ satisfy
\be \lb{qx1}
|q(x)|^{p-2} q(x)= c \sum_{i=1}^m\left(E_i\left(x;q^\pm\right)\right)^2, \qq x\in \Omega,
\ee
for some $c\ne 0$.  For later convenience, we here write the Lagrangian multiplier $c$ in the right-hand side.

For an exponent $p\in\left(1,\infty\right)$, the increasing homeomorphism $\phi_p: \R\to \R$ is given by
\[
\phi_p(s):= |s|^{p-2} s\qq \mbox{for } s\in \R,
\]
and its inverse is  $\phi_{p^*}= \phi_p^{-1}$, where $p^*:=p/(p-1) \in (1,\oo)$ is the conjugate exponent of $p$.  
\begin{lem}\label{le23}
The minimization and maximization problems of \eqref{eq66} correspond to the Lagrangian multiplier $c<0$ and $c>0$ in \eqref{qx1} respectively.
\end{lem}

\Proof
For the optimizing potentials $q=q_{m,p,r}^\pm\left(x\right)$, equation \eqref{qx1}  is equivalent to
    \[
\phi_p(q(x))=c \sum_{i=1}^m\left(E_i\left(x;q^\pm\right)\right)^2 \;\text{and}\;q(x) = \pps(c) \pps\z( \sum_{i=1}^m\left(E_i\left(x;q^\pm\right)\right)^2\y).
\]

We construct the following parameterized potentials
    \[
Q_\sigma:= \sigma\pps(c) \pps\z( \sum_{i=1}^m\left(E_i\left(x;q\right)\right)^2\y)\in \bpr,\qq\sigma\in [0,1].
\]
Clearly, $Q_1=q$. For the minimization problem \eqref{eq66}, one has
\[
\mathscr{E}_m\left(Q_\sigma\right)\geq\mathscr{E}_m\left(Q_1\right)\qqf \sigma\in[0,1].
\]
Thereby, the derivative
 \beaa
0 \GE \z. \f{\rd}{\rd\sigma} \z(\mathscr{E}_m\left(Q_\sigma\right)\y)\y|_{\sigma=1} \\
\EQ \int_\Omega\left(\sum_{i=1}^m\left(E_i\left(x;q\right)\right)^2\right)\d \pps(c)\pps\z(\sum_{i=1}^m\left(E_i\left(x;q\right)\right)^2\y) \dx\\
\EQ \pps(c) \int_\Omega \z|\sum_{i=1}^m\left(E_i\left(x;q\right)\right)^2\y|^{p^*}\dx.
\eeaa
Since $c\ne 0$, there must be $c<0$ for the minimization problem.

Analogously, $c>0$ for the maximization problem.

\qed

\noindent {\bf Proof of Theorem \ref{main}.}
 Note that $E_i(x)=E_{i}(x;q)$ are eigenfunctions for $i=1,\ldots,m$. We define the  following $m$ parameters
 $$\mu_i:= \la_i(q)\qq i=1,\ldots,m.$$
Thus,
\begin{align}\label{eq76}
-E''_i + q \left(x\right)E_i=\mu_iE_i,\; i=1,\ldots,m,\; x\in \Omega.
\end{align}

To simplify  the original critical equation \eqref{qx1}, we need to introduce the next notations
\be \lb{nw}
\e:= {\rm sign}(c) =\pm 1, \qq u_i(x):=\sqrt{|c|} E_{i}(x;q), \q i=1,\ldots,m.
\ee
Therefore, equation \eqref{qx1}  is equivalent to
\begin{align}\label{eq65}
\phi_p(q(x))=\e \sum_{i=1}^mu_i^2\;\text{and}\;q(x) = \e \pps\z( \sum_{i=1}^mu_i^2\y).
\end{align}
Since $u_i\left(x\right)$ are still eigenfunctions, by equation \eqref{eq76}, we have
$$-u''_i + q \left(x\right)u_i=\mu_iu_i, \qq i=1,\ldots,m.$$
The critical system  \eqref{ceq37} is obtained directly by substituting  \eqref{eq65} into the above system.

From the above analysis, it is not difficult to prove that equalities \eqref{u12}-\eqref{LM11} hold.  For instance, equality \eqref{qx2} is the second equality of \eqref{eq65}.  Equality \eqref{u12} is from  the norm $\nmp{q}=r$.

The proof is finished.
\qed

\section{Complete integrability}\label{se4}
In this sections,   Propositions  \ref{pr4} and \ref{pr5} show that system \eqref{eq38} is complete integrability if the parameters $k$ and $\mathscr{U}$ belong to Table \ref{Ta2}. 

%

\begin{prop}\label{pr4}
	For $\mu_1=\mu_2=\cdots=\mu_m$, the Hamiltonian system \eqref{eq38} is  completely integrable with $m-1$ additional functionally independent first integrals
	\be \lb{Integral1}
I_i=u_1v_{i+1}-u_{i+1}v_1,\quad i=1,\ldots,m-1.
	\ee
\end{prop}
\Proof
Straightforward calculations show that	$I_i=u_1v_{i+1}-u_{i+1}v_1$ are
first integrals of system \eqref{eq38} with $i=1,\ldots,m-1$. Since
$$\det\left(\partial_{\mathbf{u}}H,\partial_{\mathbf{u}}I_1,\partial_{\mathbf{u}}I_2,\ldots, \partial_{\mathbf{u}}I_{m-1}\right)=\det	\left(
\begin{array}{cccccc}
	\partial_{u_1}H&v_2&v_3&v_4&\cdots&v_m\vspace{1ex}\\
	\partial_{u_2}H&	-v_1&0&0&\cdots&0\\
	\partial_{u_2}H&	0&-v_1&0&\cdots&0\\
		\partial_{u_3}H&	0&0&-v_1&\cdots&0\\
		\partial_{u_4}H&		0&0&0&\cdots&0\\
\vdots&\vdots&\vdots&\vdots&\ddots&\vdots\\
		\partial_{u_m}H&0&0&0&\cdots&-v_1\\
\end{array}\right)\not\equiv0,$$
then $\text{rank}\left(\nabla H,\nabla I_1,\nabla I_2,\ldots,\nabla I_{m-1}\right)=m$, that is, $H$ and $I_i$
are functionally independent with $i=1,\ldots,m-1$.

\qed

For $k=2$ and $\varepsilon=+$, the Hamiltonian system \eqref{eq38} becomes a known complete integrability mechanical system, see \cite{CC77}.
\begin{lem}\label{le1}{\rm (\cite{CC77})}
 For $k=2$ and $\varepsilon=+$, the following statements hold.
 	\begin{itemize}
 	\item [\emph{(i)}] The Hamiltonian system \eqref{eq38} is  completely integrable.
 	\item[\emph{(ii)}] Let
 \begin{equation}\label{eq52}
\begin{split}
\mathcal{I}\left(\mathbf{u},\mathbf{v},\epsilon\right)=&\left(\sum_{j=1}^mu_j^2\right)\left(\sum_{j=1}^m\delta_j u_j^2\right) -\left(\sum_{j=1}^m\delta_ju_j^2\right)\left(\sum_{j=1}^m\delta_jv_j^2\right)+\\
&\left(\sum_{j=1}^m\delta_ju_jv_j\right)^2+2\sum_{j=1}^m\delta_j\left(v_j^2+\mu_ju_j^2\right)
\end{split}
 \end{equation}
with $\delta_j=1/\left(\epsilon-\mu_j\right)$  for $j=1,\ldots,m$. Then
$\mathcal{I}\left(\mathbf{u},\mathbf{v},\epsilon\right)$ is a first integral of system \eqref{eq38} independent $\epsilon$, and $\mathcal{I}\left(\mathbf{u},\mathbf{v},\epsilon\right)$ after expansion
 by powers of $\epsilon$ give $m$ functionally independent first integrals in involution for \eqref{eq39}.
\end{itemize}
\end{lem}

\begin{prop}\label{pr5}
	Let $I\left(\mathbf{u},\mathbf{v}\right)$ be a first integral of system \eqref{eq38} with $k=2$ and $\varepsilon=+$. Then the following statements hold.
	\begin{itemize}
		\item [\emph{(i)}] 	For $k=2$, the Hamiltonian system \eqref{eq38} is  completely integrable.
			\item [\emph{(ii)}] For $k=2$ and $\varepsilon=+$, $m$ functionally independent first integrals are given by statement \emph{(ii)} of Lemma \ref{le1}.
			\item [\emph{(iii)}] For $k=2$ and $\varepsilon=-$,  $I\left(-\emph{i}\mathbf{u},\emph{i}\mathbf{v}\right)$ is a first integral of system \eqref{eq38}.
		\end{itemize}
\end{prop}

\Proof
Let $I\left(\mathbf{u},\mathbf{v}\right)$ be a first integral of system \eqref{eq38} with $\varepsilon=+$.
Doing the linear canonical change of variables $$\left(\mathbf{u},\mathbf{v},t\right)\mapsto\left(\text{i}\mathbf{u},-\text{i}\mathbf{v},-t\right),$$
the integrability of the case $\varepsilon=-$ is equivalent to the case $\varepsilon=+$. Thus, $I\left(-\text{i}\mathbf{u},\text{i}\mathbf{v}\right)$ is a first integral of system \eqref{eq38} with $\varepsilon=-$. By Lemma \ref{le1}, this proposition holds.
\qed

\section{Meromorphic non-integrability}\label{se5}
In this section,  our aim is to prove the meromorphic non-integrability of Hamiltonian system \eqref{eq38} when  the parameters $k$ and $\mathscr{U}$ are outside Table \ref{Ta2}.  

\begin{prop}\label{pr7}
	If the parameters $k$ and $\mathscr{U}$ are  outside Table \ref{Ta2}, then Hamiltonian system \eqref{eq38} is meromorphic non-integrable.
\end{prop}

\Proof
Let the parameters $k$ and $\mathscr{U}$ be outside 
Table \ref{Ta2}. Then, $k\geq3$ and there exists a positive integer $j_0\in\left\{2,\ldots,m\right\}$ such that $\mu_1\neq\mu_{j_0}$. We can assume without loss of generality that $j_0=2$, that is, $\mu_1\neq\mu_2$, because in the other case  one can interchange respectively the roles of $\mu_{j_0}$ and $\mu_2$, and $u_{j_0}$ and $u_2$. To be clear, our analysis is divided into two classes:
\begin{align*}
&\text{\bf Class 1:}\;k\geq3, \mu_1\neq\mu_2\;\text{and}\; \mu_1\mu_2\neq0\;\text{(i.e. Lemma \ref{le2} below)};\\
&\text{\bf Class 2:}\;k\geq3, \mu_1\neq\mu_2\;\text{and}\; \mu_1\mu_2=0\;\text{(i.e. Lemma \ref{le3} below)}.
\end{align*}

The following Lemma \ref{le2} and Lemma \ref{le3} will complete the proof of Proposition \ref{pr7}.
\qed

\begin{lem}\label{le2}
If $k\geq3$, $\mu_1\neq\mu_2$ and $\mu_1\mu_2\neq0$, then Hamiltonian system \eqref{eq38} is meromorphic non-integrable.
\end{lem}

\Proof
One can easily observe that system \eqref{eq38} has two invariant manifolds
\begin{align*}
&\mathcal{N}_1=\left\{\left(\mathbf{u},\mathbf{v}\right)\in \mathbb{C}^{2m}\mid u_j=v_j=0,j=2,\ldots,m\right\},\\
&\mathcal{N}_2=\left\{\left(\mathbf{u},\mathbf{v}\right)\in \mathbb{C}^{2m} \mid u_1=v_1=0, u_j=v_j=0,j=3,\ldots,m\right\}.
\end{align*}
System \eqref{eq38} restricted to the first invariant manifold $\mathcal{N}_1$  becomes
\begin{align}\label{eq1}
	&u'_1=v_1,\quad v'_1=-\mu_1u_1+\varepsilon u_1^{2k-1},
\end{align}
which has first integral
\begin{align}\label{eq16}
&h=\dfrac{1}{2}v_1^2+\dfrac{1}{2}\mu_1u_1^2-\dfrac{\varepsilon}{2k}u_1^{2k}.
\end{align}
Solving equation \eqref{eq16}, we have
\begin{align}\label{eq47}
&\dfrac{du_1}{dt}=\pm\sqrt{2h+\dfrac{\varepsilon}{k}u_1^{2k}-\mu_1u_1^2}.
\end{align}
As we know, equation \eqref{eq47} for $k=2$ and $k\geq3$ is respectively called  \emph{incomplete elliptic integral of first kind} and \emph{hyperelliptic integral}, whose expressions are not always \emph{elementary functions}, see \cite{MR0060642}.

Let $\Theta\left(h\right)\in\mathbb{C}^2$ be an integral curve of system \eqref{eq1} lying on the energy level $h$.
So,
\begin{align}\label{eq57}
&\varGamma_h:=\left\{\left(u_1\left(t\right),v_1\left(t\right),0,\ldots,0\right)\in\mathbb{C}^{2m}\mid \left(u_1\left(t\right),v_1\left(t\right)\right)\in\Theta\left(h\right)\right\}
\end{align}
is a particular solution of system \eqref{eq38}.  The requirement of Theorem \ref{th5} is to construct a non-equilibrium  particular solution $\varGamma_h$. We fix the energy level  $h=0$.  Equation \eqref{eq47} has three equilibrium points   $u_1=0,\pm\sqrt[2k-2]{\mu_1k/\varepsilon}$ in the zero energy level.  To exclude these equilibrium points,  one can assume that $u_1\left(t\right)$ is not a constant. By this way, we  can get a non-equilibrium  particular solution $\Gamma_0\in\varGamma_0$.

Let
$\bm{\xi}:=\left(\xi_1,\ldots,\xi_m\right)^T$ and  $\bm{\tilde{\xi}}:=\left(\tilde{\xi}_1,\ldots,\tilde{\xi}_m\right)^T$.  We obtain that the variational equation (VE) along $\Gamma_0$ is
\begin{align}\label{eq2}
	&\left(
	\begin{array}{c}
		\bm{\xi}'\vspace{1ex}\\
		\bm{\tilde{\xi}}'\\
	\end{array}\right)
	=
	\left(
	\begin{array}{cc}
		\mathbf{0}&\mathbf{I}\vspace{2ex}\\
		\mathbf{\Lambda}&\mathbf{0}\\
	\end{array}
	\right)
	\left(
	\begin{array}{c}
		\bm{\xi}\vspace{1ex}\\
		\bm{\tilde{\xi}}\\
	\end{array}\right),
\end{align}
where
$$\mathbf{\Lambda}:=\text{diag}\;\Big(\varepsilon\left(2k-1\right)u_1^{2k-2}\left(t\right)-\mu_1,\varepsilon u_1^{2k-2}\left(t\right)-\mu_2,\varepsilon u_1^{2k-2}\left(t\right)-\mu_3,\ldots,\varepsilon u_1^{2k-2}\left(t\right)-\mu_m\Big).$$
The VE \eqref{eq2} is composed of  $m$ uncoupled Schr\"{o}dinger equations
\begin{align*}
	&\bm{\xi}''=\mathbf{\Lambda}\bm{\xi},
\end{align*}
that is,
\begin{align}
	&\xi''_1=\left(\varepsilon \left(2k-1\right)u_1^{2k-2}\left(t\right)-\mu_1\right)\xi_1,\label{eq62}\\
	&\xi''_j=\left(\varepsilon u_1^{2k-2}\left(t\right)-\mu_j\right)\xi_j,\quad j=2,\ldots,m.
\end{align}
Since $\xi_1=u'_1\left(t\right)$ is a solution of \eqref{eq62},  equation  \eqref{eq62} can be solved by Liouville's formula \cite{MR2682403}. Thereby, the normal variational equations (NVE) along $\Gamma_0$ are
\begin{align}\label{eq3}
	&\xi''_j=\left(\varepsilon u_1^{2k-2}\left(t\right)-\mu_j\right)\xi_j,\quad j=2,\ldots,m.
\end{align}

Inspired
by  Yoshida \cite{MR923886}, we introduce the following finite branched covering map
\begin{equation}\label{eq17}
	\begin{split}
	&\overline{\Gamma}_0\rightarrow \mathbf{P}^1,\\
	&t\longmapsto z=\dfrac{\varepsilon}{k\mu_1}u_1^{2k-2}\left(t\right),
	\end{split}
\end{equation}
where $\overline{\Gamma}_0$ is the compact Riemann surface of the curve $v_1^2=\varepsilon u_1^{2k}/k-\mu_1u_1^2$ and $\mathbf{P}^1$ is the Riemann sphere.

Performing the Yoshida transformation \eqref{eq17}, the normal variational equations \eqref{eq3} can be written as the hypergeometric differential equations in the new independent variable $z$
\begin{align}
	&\dfrac{d^2\xi_j}{dz^2}+\left(\frac{1}{z}+\frac{1}{2 (z-1)}\right)\dfrac{d\xi_j}{dz}-\left(\frac{\mu _j}{4 \mu _1(k-1)^2 z^2}+\frac{k \mu _1-\mu _j}{4\mu _1 (k-1)^2z (z-1)}\right)\xi_j=0,\tag{$\text{ANVE}_j$}\\
	&j=2,\ldots,m.\notag
\end{align}
The above differential system of equations is called the \emph{algebraic normal variational equations} (ANVE), and is denoted as
\begin{align}\label{eq4}
	&\text{ANVE}=\text{ANVE}_2\oplus\text{ANVE}_3\oplus\cdots\oplus\text{ANVE}_m.
\end{align}
Essentially, equation \eqref{eq4} is a direct sum in the more intrinsic sense of linear connections, see Chapter $2$ of \cite{MR1713573}  for more details.

From Theorem \ref{th4}, it follows that the identity components of the Galois groups of the NVE \eqref{eq3} and the ANVE \eqref{eq4} coincide. Obviously, the ANVE \eqref{eq4}  is integrable if and only if each $\text{ANVE}_j$ is integrable for $j=1,\ldots,m$. More precisely,  the identity component of the
Galois group of the ANVE is solvable if and only if the identity
 component of the Galois group of each $\text{ANVE}_j$ is solvable  for $j=1,\ldots,m$.

Now, we consider the $\text{ANVE}_2$:
\begin{align}\label{eq19}
&\dfrac{d^2\xi_2}{dz^2}+\left(\frac{1}{z}+\frac{1}{2 (z-1)}\right)\dfrac{d\xi_2}{dz}-\left(\frac{\mu _2}{4 \mu _1(k-1)^2 z^2}+\frac{k \mu _1-\mu _2}{4\mu _1 (k-1)^2z (z-1)}\right)\xi_2=0
\end{align}
with three singular points at $z=0,1,\infty$. Comparing  \eqref{eq19} with the general form of the
hypergeometric equation \eqref{hy},  one can see that the exponents of \eqref{eq19} at singular points must
fulfill the following relations
\begin{align*}
&\alpha+\tilde{\alpha}=0,\quad \alpha\tilde{\alpha}=-\frac{\mu _2}{4 \mu _1(k-1)^2},\\
&\beta+\tilde{\beta}=\dfrac{1}{2},\quad \beta\tilde{\beta}=-\frac{k}{4(k-1)^2},\\
&\gamma+\tilde{\gamma}=\dfrac{1}{2},\quad \gamma\tilde{\gamma}=0.
\end{align*}
Thus, all the possibilities of the differences of exponents are
\begin{align}\label{eq21}
&\varrho=\pm\dfrac{1}{k-1}\sqrt{\dfrac{\mu_2}{\mu_1}},\tau=\pm\dfrac{1}{2}\left(1+\dfrac{2}{k-1}\right)\;\text{and}\;\varsigma=\pm\dfrac{1}{2}.
\end{align}
Moreover, we can get all the possibilities of $\varrho+\tau+\varsigma,-\varrho+\tau+\varsigma,\varrho-\tau+\varsigma$ and $\varrho+\tau-\varsigma$, see Table \ref{Ta1}.
\begin{table}[H]
\tiny
	\centering
	\caption{All the possibilities of $\varrho+\tau+\varsigma,-\varrho+\tau+\varsigma,\varrho-\tau+\varsigma$ and $\varrho+\tau-\varsigma$.}\vspace{2ex}
	\begin{tabular}{|c|c|c|c|c|c|}
		\hline
		&Signs of $\left(\varrho,\tau,\varsigma\right)$&$\varrho+\tau+\varsigma$&$-\varrho+\tau+\varsigma$&	$\varrho-\tau+\varsigma$&$\varrho+\tau-\varsigma$\xrowht{20pt}\\
		\hline
		$1$&$\left(+,+,+\right)$&$1+\dfrac{1}{k-1}\left(\sqrt{\dfrac{\mu_2}{\mu_1}}+1\right)$&$1+\dfrac{1}{k-1}\left(1-\sqrt{\dfrac{\mu_2}{\mu_1}}\right)$&$\dfrac{1}{k-1}\left(\sqrt{\dfrac{\mu_2}{\mu_1}}-1\right)$&$\dfrac{1}{k-1}\left(\sqrt{\dfrac{\mu_2}{\mu_1}}+1\right)$\xrowht{28pt}\\
		\hline
		$2$&$\left(+,+,-\right)$&$\dfrac{1}{k-1}\left(\sqrt{\dfrac{\mu_2}{\mu_1}}+1\right)$&$\dfrac{1}{k-1}\left(1-\sqrt{\dfrac{\mu_2}{\mu_1}}\right)$&$\dfrac{1}{k-1}\left(\sqrt{\dfrac{\mu_2}{\mu_1}}-1\right)-1$&$1+\dfrac{1}{k-1}\left(\sqrt{\dfrac{\mu_2}{\mu_1}}+1\right)$\xrowht{28pt}\\
		\hline
		$3$&$\left(+,-,+\right)$&$\dfrac{1}{k-1}\left(\sqrt{\dfrac{\mu_2}{\mu_1}}-1\right)$&$-\dfrac{1}{k-1}\left(\sqrt{\dfrac{\mu_2}{\mu_1}}+1\right)$&$1+\dfrac{1}{k-1}\left(\sqrt{\dfrac{\mu_2}{\mu_1}}+1\right)$&$\dfrac{1}{k-1}\left(\sqrt{\dfrac{\mu_2}{\mu_1}}-1\right)-1$\xrowht{28pt}\\
		\hline
		$4$&$\left(+,-,-\right)$&$\dfrac{1}{k-1}\left(\sqrt{\dfrac{\mu_2}{\mu_1}}-1\right)-1$&$-1-\dfrac{1}{k-1}\left(\sqrt{\dfrac{\mu_2}{\mu_1}}+1\right)$&$\dfrac{1}{k-1}\left(\sqrt{\dfrac{\mu_2}{\mu_1}}+1\right)$&$\dfrac{1}{k-1}\left(\sqrt{\dfrac{\mu_2}{\mu_1}}-1\right)$\xrowht{28pt}\\
		\hline
		$5$&$\left(-,+,+\right)$&$1+\dfrac{1}{k-1}\left(1-\sqrt{\dfrac{\mu_2}{\mu_1}}\right)$&$1+\dfrac{1}{k-1}\left(\sqrt{\dfrac{\mu_2}{\mu_1}}+1\right)$&$-\dfrac{1}{k-1}\left(\sqrt{\dfrac{\mu_2}{\mu_1}}+1\right)$&$\dfrac{1}{k-1}\left(1-\sqrt{\dfrac{\mu_2}{\mu_1}}\right)$\xrowht{28pt}\\
		\hline
		$6$&$\left(-,+,-\right)$&$\dfrac{1}{k-1}\left(1-\sqrt{\dfrac{\mu_2}{\mu_1}}\right)$&$\dfrac{1}{k-1}\left(\sqrt{\dfrac{\mu_2}{\mu_1}}+1\right)$&$-1-\dfrac{1}{k-1}\left(\sqrt{\dfrac{\mu_2}{\mu_1}}+1\right)$&1+$\dfrac{1}{k-1}\left(1-\sqrt{\dfrac{\mu_2}{\mu_1}}\right)$\xrowht{28pt}\\
		\hline
		$7$&$\left(-,-,+\right)$&$-\dfrac{1}{k-1}\left(\sqrt{\dfrac{\mu_2}{\mu_1}}+1\right)$&$\dfrac{1}{k-1}\left(\sqrt{\dfrac{\mu_2}{\mu_1}}-1\right)$&$1+\dfrac{1}{k-1}\left(1-\sqrt{\dfrac{\mu_2}{\mu_1}}\right)$&$-1-\dfrac{1}{k-1}\left(\sqrt{\dfrac{\mu_2}{\mu_1}}+1\right)$\xrowht{28pt}\\
		\hline
		$8$&$\left(-,-,-\right)$&$-1-\dfrac{1}{k-1}\left(\sqrt{\dfrac{\mu_2}{\mu_1}}+1\right)$&$\dfrac{1}{k-1}\left(\sqrt{\dfrac{\mu_2}{\mu_1}}-1\right)-1$&$\dfrac{1}{k-1}\left(1-\sqrt{\dfrac{\mu_2}{\mu_1}}\right)$&$-\dfrac{1}{k-1}\left(\sqrt{\dfrac{\mu_2}{\mu_1}}+1\right)$\xrowht{28pt}\\
		\hline
	\end{tabular}%
	\label{Ta1}%
\end{table}%

If equation \eqref{eq19} satisfies statement (i) of  Theorem \ref{A1}, by Table \ref{Ta1}, then
$$\dfrac{1}{k-1}\left(\sqrt{\dfrac{\mu_2}{\mu_1}}+1\right)\;\text{or}\;\dfrac{1}{k-1}\left(\sqrt{\dfrac{\mu_2}{\mu_1}}-1\right)$$
must be an integer, that is,
\begin{align}\label{eq20}
&\dfrac{\mu_2}{\mu_1}\in\left\{\left(\left(k-1\right)\ell\pm1\right)^2\big|\ell\in \mathbb{N}\right\}.
\end{align}

The statement (ii) of  Theorem \ref{A1} has $15$ possibilities in the Table \ref{Ta3}.
If the statement (ii) of  Theorem \ref{A1} is fulfilled for equation \eqref{eq19}, from equation \eqref{eq21}, we find that only the first row of Table \ref{Ta3} conforms.
Note that $k\geq3$. Therefore,
$$\pm\dfrac{1}{k-1}\sqrt{\dfrac{\mu_2}{\mu_1}}=\dfrac{1}{2}+\ell,\quad \ell\in\mathbb{Z},$$
that is,
$$\dfrac{\mu_2}{\mu_1}\in\left\{\dfrac{\left(k-1\right)^2\left(2\ell+1\right)^2}{4}\Bigg|\ell\in \mathbb{Z}\right\}.$$

Based on the analysis above,  the parameters $\mu_1$ and $\mu_2$ must satisfy
\begin{align}\label{eq22}
&\dfrac{\mu_2}{\mu_1}\in\left\{\left(\left(k-1\right)\ell\pm1\right)^2\big|\ell\in \mathbb{N}\right\}\bigcup\left\{\dfrac{\left(k-1\right)^2\left(2\ell+1\right)^2}{4}\Bigg|\ell\in \mathbb{Z}\right\}
\end{align}
if the identity components of the Galois groups of the NVE \eqref{eq3} is Abelian.

On the second invariant manifold $\mathcal{N}_2$, system \eqref{eq38} is written as
\begin{align}\label{eq6}
	&u'_2=v_2,\quad v'_2=-\mu_2u_2+\varepsilon u_2^{2k-1}
\end{align}
with Hamiltonian
\begin{align}\label{eq23}
&\tilde{h}=\dfrac{1}{2}v_2^2+\dfrac{1}{2}\mu_2u_2^2-\dfrac{\varepsilon}{2k}u_2^{2k}.
\end{align}
To solve equation \eqref{eq23}, we
\begin{align}\label{eq48}
	&\dfrac{du_2}{dt}=\pm\sqrt{2\tilde{h}+\dfrac{\varepsilon}{k}u_2^{2k}-\mu_2u_2^2}.
\end{align}

Let $\widetilde{\Theta}\left(\tilde{h}\right)\in\mathbb{C}^2$ be a integral curve of system \eqref{eq6} lying on the energy level $\tilde{h}$. Thus,
\begin{align}\label{eq58}
	&\widetilde{\varGamma}_{\tilde{h}}:=\left\{\left(0,0,u_2\left(t\right),v_2\left(t\right),0,\ldots,0\right)\in\mathbb{C}^{2m}\mid \left(u_2\left(t\right),v_2\left(t\right)\right)\in\widetilde{\Theta}\left(\tilde{h}\right)\right\}
\end{align}
is a particular solution of system \eqref{eq38}. We select the energy level  $\tilde{h}=0$.  Equation \eqref{eq48} has three equilibrium points   $u_2=0,\pm\sqrt[2k-2]{\mu_2k/\varepsilon}$ in the zero energy level. In the same way as particular solution $\Gamma_0$, we can find a non-equilibrium  particular solution $\widetilde{\Gamma}_0\in\widetilde{\varGamma}_0$.

Let $\bm{\eta}:=\left(\eta_1,\ldots,\eta_m\right)^T$ and $\bm{\tilde{\eta}}:=\left(\tilde{\eta}_1,\ldots,\tilde{\eta}_m\right)^T$. The variational equations (VE) along $\widetilde{\Gamma}_0$ is given by
\begin{align}\label{eq7}
	&\left(
	\begin{array}{c}
		\bm{\eta}'\vspace{1ex}\\
		\bm{\tilde{\eta}}'\\
	\end{array}\right)
	=
	\left(
	\begin{array}{cc}
		\mathbf{0}&\mathbf{I}\vspace{2ex}\\
		\mathbf{\tilde{\Lambda}}&\mathbf{0}\\
	\end{array}
	\right)
	\left(
	\begin{array}{c}
		\bm{\eta}\vspace{1ex}\\
	\bm{\tilde{\eta}}\\
	\end{array}\right),
\end{align}
where
$$\mathbf{\tilde{\Lambda}}:=\text{diag}\;\Big(\varepsilon u_2^{2k-2}\left(t\right)-\mu_1,\varepsilon\left(2k-1\right)u_2^{2k-2}\left(t\right)-\mu_2,\varepsilon u_2^{2k-2}\left(t\right)-\mu_3,\ldots, \varepsilon u_2^{2k-2}\left(t\right)-\mu_m\Big).$$
The VE \eqref{eq7} is also composed of  $m$ uncoupled Schr\"{o}dinger equations
\begin{align*}
	&\bm{\eta}''=\mathbf{\tilde{\Lambda}}\bm{\eta},
\end{align*}
that is,
\begin{align}
	&\eta''_1=\left(\varepsilon u_2^{2k-2}\left(t\right)-\mu_1\right)\eta_1,\notag\\
	&\eta''_2=\left(\varepsilon\left(2k-1\right)u_2^{2k-2}\left(t\right)-\mu_2\right)\eta_2,\label{eq64}\\
	&\eta''_j=\left(\varepsilon u_2^{2k-2}\left(t\right)-\mu_j\right)\eta_j,\quad j=3,\ldots,m.\notag
\end{align}
Using Liouville's formula \cite{MR2682403}, the second equation of \eqref{eq64} is solvable due to the fact that it has a solution $\eta_2=u'_2\left(t\right)$. Therefore, the corresponding normal variational equations ($\widetilde{\text{NVE}}$) along $\widetilde{\Gamma}_0$ are given by
\begin{equation}\label{eq8}
	\begin{split}
	&\eta''_1=\left(\varepsilon u_2^{2k-2}\left(t\right)-\mu_1\right)\eta_1,\\
&\eta''_j=\left(\varepsilon u_2^{2k-2}\left(t\right)-\mu_j\right)\eta_j,\quad j=3,\ldots,m.
	\end{split}
\end{equation}

Similarly, we can carry out the following Yoshida transformation
$$t\longmapsto z=\dfrac{\varepsilon}{k\mu_2}u_2^{2k-2}\left(t\right),$$
and transform $\widetilde{\text{NVE}}$ \eqref{eq8} into the algebraic normal variational equations ($\widetilde{\text{ANVE}}$):
\begin{align}
	&\dfrac{d^2\eta_1}{dz^2}+\left(\frac{1}{z}+\frac{1
	}{2 (z-1)}\right)\dfrac{d\eta_1}{dz}-\left(\frac{\mu _1}{4 \mu _2(k-1)^2 z^2}+\frac{k \mu _2-\mu _1}{4\mu _2 (k-1)^2z (z-1)}\right)\eta_1=0,\tag{$\widetilde{\text{ANVE}}_1$}\\
	&\dfrac{d^2\eta_j}{dz^2}+\left(\frac{1}{z}+\frac{1}{2 (z-1)}\right)\dfrac{d\eta_j}{dz}-\left(\frac{\mu _j}{4 \mu _2(k-1)^2 z^2}+\frac{k \mu _2-\mu _j}{4\mu _2 (k-1)^2z (z-1)}\right)\eta_j=0,\tag{$\widetilde{\text{ANVE}}_j$}\\
	&j=3,\ldots,m.\notag
\end{align}
The direct sum form of $\widetilde{\text{ANVE}}$ is
$\widetilde{\text{ANVE}}=\widetilde{\text{ANVE}}_1\oplus\widetilde{\text{ANVE}}_3\oplus\widetilde{\text{ANVE}}_4\oplus\cdots\oplus\widetilde{\text{ANVE}}_m$. For the $\widetilde{\text{ANVE}}_1$,  all the possibilities of the differences of exponents are
$$\varrho=\pm\dfrac{1}{k-1}\sqrt{\dfrac{\mu_1}{\mu_2}},\tau=\pm\dfrac{1}{2}\left(1+\dfrac{2}{k-1}\right)\;\text{and}\;\varsigma=\pm\dfrac{1}{2}.$$

 By the same discussions as NVE \eqref{eq3}, we obtain  that  the parameters $\mu_1$ and $\mu_2$ must satisfy
 \begin{align}\label{eq24}
& \dfrac{\mu_1}{\mu_2}\in\left\{\left(\left(k-1\right)\ell\pm1\right)^2\big|\ell\in \mathbb{N}\right\}\bigcup\left\{\dfrac{\left(k-1\right)^2\left(2\ell+1\right)^2}{4}\Bigg|\ell\in \mathbb{Z}\right\} 	
 \end{align}
if the identity components of the Galois groups of the $\widetilde{\text{NVE}}$ \eqref{eq8} is Abelian.

The conditions \eqref{eq22} and  \eqref{eq24} imply that
$$\dfrac{\mu_2}{\mu_1}\geq1\;\text{and}\;\dfrac{\mu_1}{\mu_2}\geq1,$$
respectively. This
contradicts our assumption $\mu_1\neq\mu_2$. Consequently, either the identity components of the Galois groups of the NVE \eqref{eq3} or  $\widetilde{\text{NVE}}$ \eqref{eq8} is not Abelian.
By Theorem \ref{th5}, the Hamiltonian system \eqref{eq38}  for $k\geq3$  is meromorphic non-integrable with $\mu_1\neq\mu_2$ and $\mu_1\mu_2\neq0$.

The proof is finished.
\qed

\begin{lem}\label{le3}
If $k\geq3$, $\mu_1\neq\mu_2$ and $\mu_1\mu_2=0$, then Hamiltonian system \eqref{eq38} is meromorphic non-integrable.
\end{lem}

\Proof Our proof will be distinguished two cases:
$$\text{\bf Case 1:}\; \mu_1=0, \mu_2\neq0\;\text{and}\;\text{\bf Case 2:}\; \mu_1\neq0, \mu_2=0.$$

  {\bf Case 1: $\mu_1=0$ and $\mu_2\neq0$.} For this case, we also restrict system \eqref{eq38} on the invariant manifold $\mathcal{N}_1$. Namely,
\begin{align}\label{eq15}
	&u'_1=v_1,\quad v'_1=\varepsilon u_1^{2k-1}
\end{align}
with Hamiltonian
\begin{align}\label{eq25}
	&h=\dfrac{1}{2}v_1^2-\dfrac{\varepsilon}{2k}u_1^{2k}.
\end{align}

Analogously, we also consider the particular solution $\Gamma_0$ in the proof of Lemma \ref{le2}, and compute the $\widehat{\text{NVE}}$ along $\Gamma_0$
\begin{equation}\label{eq45}
	\begin{split}
&\xi''_j=\left(\varepsilon u_1^{2k-2}\left(t\right)-\mu_j\right)\xi_j,\quad j=2,\ldots,m.
	\end{split}
\end{equation}

Doing the change of variable
$$t\longmapsto z=\dfrac{\varepsilon}{2\mu_2}u_1^{2k-2}\left(t\right),$$
we attain the algebraic normal variational equations ($\widehat{\text{ANVE}}$):
\begin{align}
	&\dfrac{d^2\xi_2}{d z^2}+\dfrac{3}{2z}\dfrac{d\xi_2}{d z}-\frac{k\left(2z-1\right)}{8 (k-1)^2z^3}\xi_2=0,\tag{$\widehat{\text{ANVE}}_2$}\\
	&\dfrac{d^2\xi_j}{d z^2}+\dfrac{3}{2z}\dfrac{d\xi_j}{d z}-\frac{k \left(2\mu _2z-\mu _j\right)}{8(k-1)^2\mu _2 z^3}\xi_j=0,\quad j=3,\ldots, m, \tag{$\widehat{\text{ANVE}}_j$}
\end{align}
and denote by
\begin{align}\label{eq5}
	&\widehat{\text{ANVE}}=\widehat{\text{ANVE}}_2\oplus\widehat{\text{ANVE}}_3\oplus\cdots\oplus\widehat{\text{ANVE}}_m.
\end{align}
Making the classical transformation
$$\xi_2=\chi\exp\left(-\dfrac{3}{4}\int \dfrac{dz}{z} \right)=\chi z^{-3/4},$$
the $\widehat{\text{ANVE}}_2$ reads
\begin{align}\label{eq26}
&\chi''=r\left(z\right)\chi,
\end{align}
where
\begin{align}\label{eq27}
&r\left(z\right)=-\left(\dfrac{(k-3) (3 k-1)}{16 (k-1)^2 z^2}+\dfrac{ k}{8(k-1)^2 z^3}\right).
\end{align}

 Then, the set of poles of  $r\left(z\right)$ is $\Upsilon=\left\{0,\infty\right\}$.  The order of  $z=0$ and $z=\infty$ is $o\left(0\right)=3$ and $o\left(\infty\right)=2$, respectively.  Using Proposition \ref{pr1} to equation  \eqref{eq26}, only types (ii) or (iv) of Theorem \ref{th6} can appear. Working the second part of Kovacic's algorithm (see Appendix \ref{B}), we obtain that
$$\mathcal{E}_\infty=\left\{2+\ell\sqrt{1-\dfrac{(k-3) (3 k-1)}{4(k-1)^2 }}\;\Bigg|\;\ell=0,\pm2\right\}\bigcap\mathbb{Z}=
\begin{cases}
	\left\{0,2,4\right\},\;\text{if}\;k=3,\\
	\left\{2\right\},\;\text{if}\;k\geq4.
\end{cases}\;\text{and}\;\mathcal{E}_0=\left\{3\right\}.$$
Straightforward computations show that the  number $d=d\left(\bm{\varpi}\right)=\left(\varpi_\infty-\sum_{c\in\Upsilon}\varpi_c\right)/2$
is not a non-negative integer. Therefore, type  (iv) of Theorem \ref{th6}  holds. This means that the identity component of the
Galois group of the $\widehat{\text{ANVE}}$ \eqref{eq5} is not Abelian. Thereby, the identity component of the
Galois group of the  $\widehat{\text{NVE}}$ \eqref{eq45} is also not Abelian. From Theorem \ref{th5}, it follows that the Hamiltonian system \eqref{eq38} for $k\geq3$ is meromorphic non-integrable with $\mu_1=0$ and $\mu_2\neq0$.

 {\bf Case 2: $\mu_1\neq0$ and $\mu_2=0$.} Substituting $\mu_2=0$ into \eqref{eq8}, we get the normal variational equations along $\widetilde{\Gamma}_0$:
 \begin{equation}\label{eq46}
 	\begin{split}
&\eta''_1=\left(\varepsilon u_2^{2k-2}\left(t\right)-\mu_1\right)\eta_1,\\
&\eta''_j=\left(\varepsilon u_2^{2k-2}\left(t\right)-\mu_j\right)\eta_j,\quad j=3,\ldots,m.
 	\end{split}
 \end{equation}
After the change of variable
$$t\longmapsto z=\dfrac{\varepsilon}{2\mu_1}u_2^{2k-2}\left(t\right),$$
equations \eqref{eq46} become the algebraic normal variational equations
 \begin{align*}
 	&\dfrac{d^2\eta_1}{d z^2}+\dfrac{3}{2z}\dfrac{d\eta_1}{d z}-\frac{k\left(2z-1\right)}{8 (k-1)^2z^3}\eta_1=0,\\
 	&\dfrac{d^2\eta_j}{d z^2}+\dfrac{3}{2z}\dfrac{d\eta_j}{d z}-\frac{k \left(2\mu _1z-\mu _j\right)}{8(k-1)^2\mu _1 z^3}\eta_j=0,\quad j=3,\ldots, m.
 \end{align*}
  The analysis is exactly the same as {\bf Case 1}.  Thus, the Hamiltonian system \eqref{eq38} for $k\geq3$ is meromorphic non-integrable with $\mu_1\neq0$ and $\mu_2=0$.

 This lemma holds.
\qed

\noindent {\bf Proof of Theorem \ref{th1}.}
 By Propositions \ref{pr4}, \ref{pr5} and \ref{pr7}, the Theorem \ref{th1} follows.
\qed

\section{Poincar\'{e} cross section}\label{se6}
For  integrable Hamiltonian systems, Liouville-Arnold theorem \cite{MR2004534} exhibits that their dynamical  behaviors  are ordered and regular.  For weakly perturbed (originally integrable) Hamiltonian systems, KAM  theorem \cite{MR2004534} shows their dynamical  behaviors  are stochastic and chaotic, such as chaos, Arnold diffusion (at least three degrees of freedom), etc.
Roughly speaking, the phase space for
integrable Hamiltonian systems
is foliated by KAM tori, which obstruct the stochasticity of the trajectories. For some weak perturbations, the KAM tori happen breaking down  resulting in chaotic motion. In other words, the chaotic behavior can destroy meromorphic integrability.  The classical Poincar\'{e} cross section technique can intuitively present the above dynamical  process: local stability, the trajectories transition from ordered to chaotic, and many other dynamic properties.

\begin{figure}[H]
	\centering
	\begin{minipage}{0.45\linewidth}
		\centering
		\centerline{\includegraphics[width=1.1\textwidth]{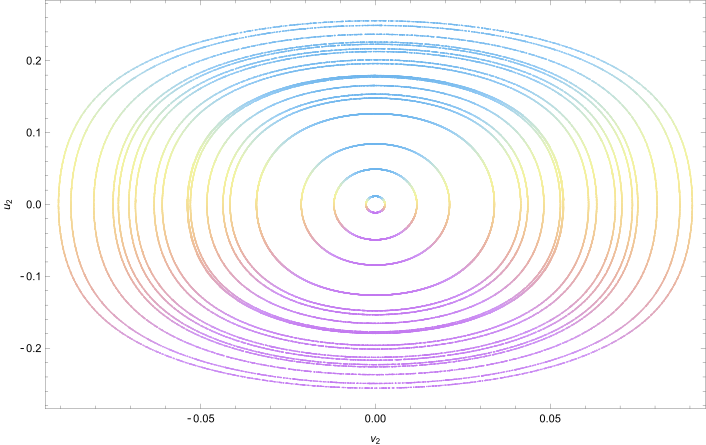}}
		\centerline{(1) $\mu_2=1$.}
	\end{minipage}
	\quad
	\quad
	\begin{minipage}{0.45\linewidth}
		\centering
		\centerline{\includegraphics[width=1.1\textwidth]{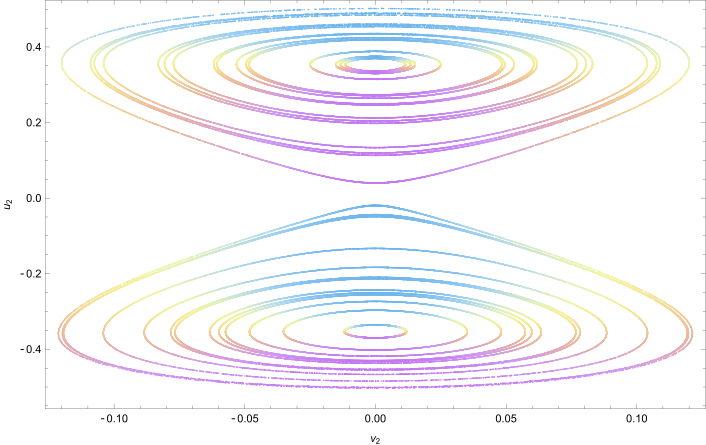}}
		\centerline{(2) $\mu_2=-1$.}
	\end{minipage}
	\caption{Poincar\'{e} cross section $\left(v_2,u_2\right)$: the integrable case $k=2$, $\mu_1=0.1$ and $h=0.85$.}\label{fig2}
\end{figure}

In the calculation Poincar\'{e} cross sections below, we focus on Hamiltonian system \eqref{eq38}  with two degrees of freedom (i.e. $m=2$) and fix the  parameter $\varepsilon=-1$. Consider the energy level
$$M_h:=\left\{\left(u_1,u_2,v_1,v_2\right)\in\mathbb{R}^4| H\left(u_1,u_2,v_1,v_2\right)=h,\;h\in\mathbb{R}\right\}.$$
On energy level $M_h$, we  select $u_1=0$ as a
cross section plane with coordinates $\left(u_2,v_2\right)$.
Taking the  energy $h=0.85$, Figure \ref{fig2} shows the integrable case with $k=2$, which take the set of parameters: $\mu_1=0.1$, $\mu_2=1$; $\mu_1=0.1$, $\mu_2=-1$.  We can see that their dynamical structures are very regular.

\begin{figure}[H]
	\centering
	\begin{minipage}{0.45\linewidth}
		\centering
		\centerline{\includegraphics[width=1.1\textwidth]{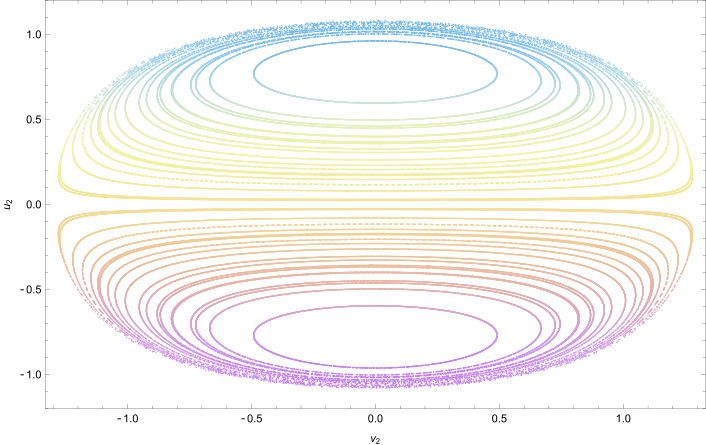}}
		\centerline{(1) $\mu_1=1$.}
	\end{minipage}
	\quad
	\quad
	\begin{minipage}{0.45\linewidth}
		\centering
		\centerline{\includegraphics[width=1.1\textwidth]{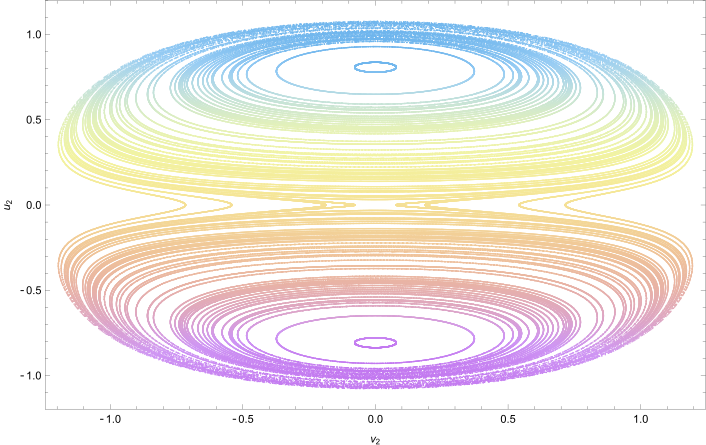}}
		\centerline{(2) $\mu_1=0.99$.}
	\end{minipage}

	\begin{minipage}{0.45\linewidth}
	\centering
	\centerline{\includegraphics[width=1.1\textwidth]{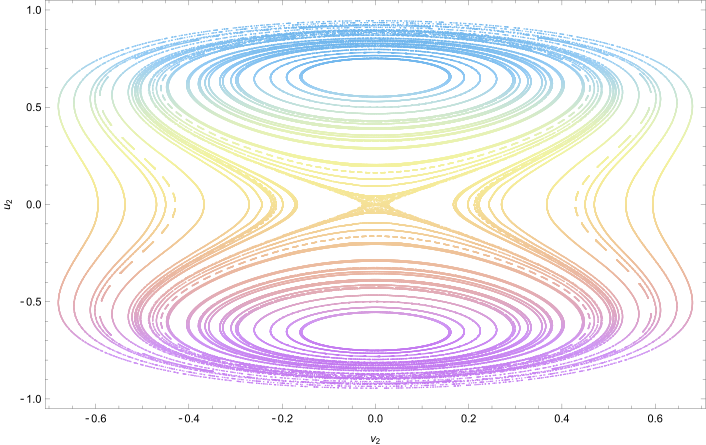}}
	\centerline{(3) $\mu_1=0.5$.}
\end{minipage}
	\quad
\quad
	\begin{minipage}{0.45\linewidth}
		\centering
		\centerline{\includegraphics[width=1.1\textwidth]{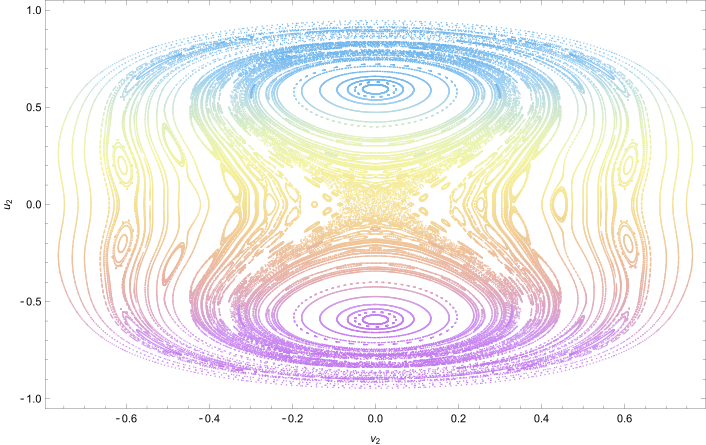}}
		\centerline{(4) $\mu_1=0.1$.}
	\end{minipage}
	\caption{Poincar\'{e} cross section $\left(v_2,u_2\right)$:
$k=3$, $\mu_2=1$ and $h=0.85$.}\label{fig1}
\end{figure}

The integrable case  and some weakly perturbed cases  with $k=3$, $\mu_2=1$ and $h=0.85$ are presented in Figure \ref{fig1}. For integrable case,  the dynamical behavior is highly regular, see (1) of Figure \ref{fig1}.  For sufficiently small perturbation, KAM tori appear deformation or even breaking down in the fragile top and bottom boundaries, but most of them remain for internal region, as shown in (2) and (3) of Figure \ref{fig1}. For $\mu_1=0.99$ and $\mu_1=0.5$, the major bifurcations occurs in the vertical direction.  As the perturbation strength is increased $\mu_1=0.1$, the KAM tori progressively break down resulting in the trajectories complete stochastic motion.  Accurately speaking, the central domain in (4) of Figure \ref{fig1} is a  large chaotic
zone, that is, chaos. Around this chaotic zone, there are many chain of islands  which
correspond to quasi-periodic trajectories. Except for periodic trajectories, some KAM tori still remain in the annular area at top and bottom, which can be observed in (4) of Figure \ref{fig1}.

\begin{figure}[H]
	\centering
\begin{minipage}{0.45\linewidth}
	\centering
	\centerline{\includegraphics[width=1.1\textwidth]{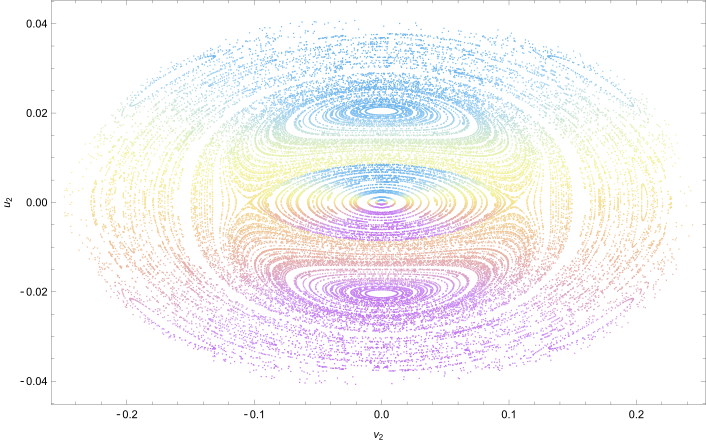}}
	\centerline{(1) $k=10$, $\mu_2=45$.}
\end{minipage}
\quad
\quad
\begin{minipage}{0.45\linewidth}
	\centering
	\centerline{\includegraphics[width=1.1\textwidth]{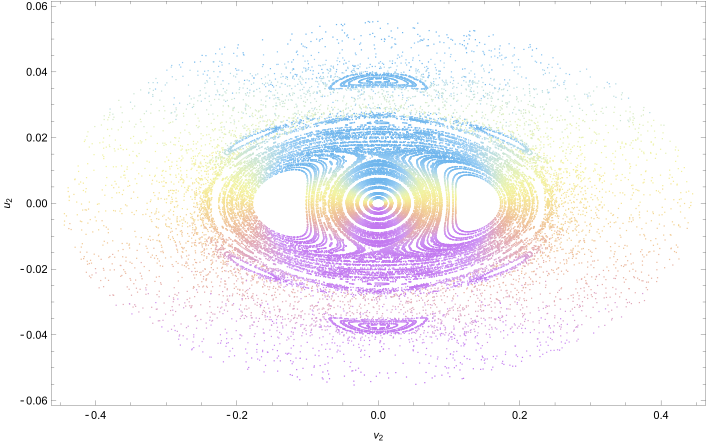}}
	\centerline{(2) $k=20$, $\mu_2=65$.}
\end{minipage}

\begin{minipage}{0.45\linewidth}
	\centering
	\centerline{\includegraphics[width=1.1\textwidth]{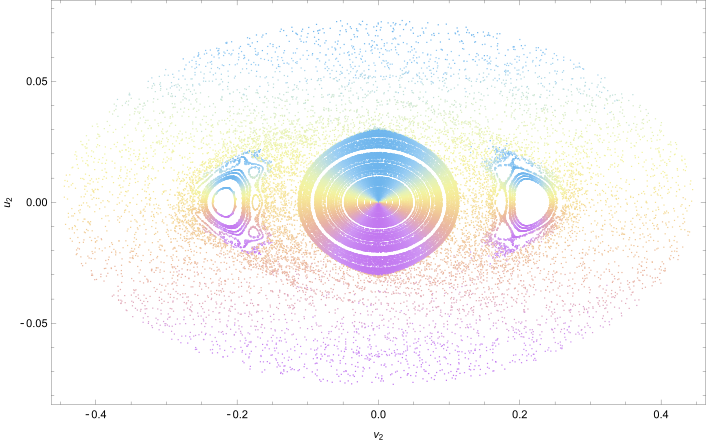}}
	\centerline{(3) $k=30$, $\mu_2=35$.}
\end{minipage}
\quad
\quad
\begin{minipage}{0.45\linewidth}
	\centering
	\centerline{\includegraphics[width=1.1\textwidth]{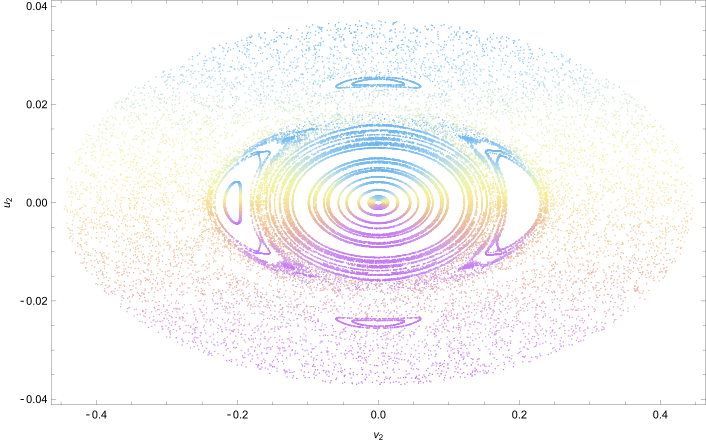}}
	\centerline{(4) $k=40$, $\mu_2=115$.}
\end{minipage}
	\caption{Poincar\'{e} cross section $\left(v_2,u_2\right)$: for some high degree $k$ with $\mu_1=-1.5$ and $h=0.1$.}\label{fig3}
\end{figure}

Finally, our numerical simulations consider some  high degree systems \eqref{eq38}, that is, $k=10,20,30,40$, see Figure \ref{fig3}.  One can see that the KAM tori  at central region boundary disappear with the trajectories escaping to nonclosed areas of the phase space.   This trajectories stochastic escaping create  complex dynamic phenomena of system \eqref{eq38},  including chaos, quasi-periodic trajectories and periodic trajectories.

\section*{Acknowledgments}
Yuzhou Tian wants to express his gratitude to the Department of Mathematical Sciences, Tsinghua University for the hospitality and support during the time period in which this work was completed.

The second author is partially supported by the National Natural Science Foundation of China (Grant no. 11790273).

\section*{Appendix}
\setcounter{equation}{0}
\setcounter{subsection}{0}
\setcounter{figure}{0}
\setcounter{table}{0}
\renewcommand{\theequation}{\thesection.\arabic{equation}}
\renewcommand{\thesubsection}{\thesection.\arabic{subsection}}
\renewcommand{\thetable}{\thesection.\arabic{table}}
\renewcommand{\thefigure}{\thesection.\arabic{figure}}
\begin{appendices}
\section{Second part of Kovacic’s algorithm}\label{B}
Here, we recall  the second part of  Kovacic's algorithm \cite{MR839134}. Let $r\left(z\right)\in C\left(z\right)$ and $\Upsilon$ be the set of poles of $r\left(z\right)$. Set $\chi''=r\chi$.
	\begin{itemize}
		\item [\bf{Step 1.}] To each pole $c\in\Upsilon$, we calculate the set  $\mathcal{E}_c$ as follows.
		\begin{itemize}
			\item [(i)] If the pole $c$ is of order $1$, then $\mathcal{E}_c=\left\{4\right\}$.
			\item [(ii)] If the pole $c$ is of order $2$ and $b$ is  the coefficient of $1/\left(z-c\right)^2$ in the partial fraction decomposition of $r\left(z\right)$, then
			\begin{align}\label{eq29}
				&\mathcal{E}_c=\left\{2+\ell\sqrt{1+4b}\;\big|\;\ell=0,\pm2\right\}\bigcap\mathbb{Z}.
			\end{align}
			\item [(iii)] If the pole $c$ is of order $o\left(c\right)>2$, then $\mathcal{E}_c=\left\{o\left(c\right)\right\}$.
			\item [(iv)] If the order of $r$ at $\infty$ is $o\left(\infty\right)>2$, then $\mathcal{E}_c=\left\{0,2,4\right\}$.
			\item [(v)] If the order of $r$ at $\infty$ is $2$ and $b$ is  the coefficient of $1/z^2$ in the Laurent expansion  of $r\left(z\right)$ at $\infty$, then
			\begin{align}\label{eq42}
				&\mathcal{E}_c=\left\{2+\ell\sqrt{1+4b}\;\big|\;\ell=0,\pm2\right\}\bigcap\mathbb{Z}.
			\end{align}
			\item [(vi)] If the order of $r$ at $\infty$ is $o\left(\infty\right)<2$, then $\mathcal{E}_c=\left\{o\left(\infty\right)\right\}$.
		\end{itemize}
		\item [\bf{Step 2.}] Let  $\bm{\varpi}=\left(\varpi_c\right)_{c\in\Upsilon}$ be a element in the Cartesian product $\prod_{c\in\Upsilon}\mathcal{E}_c$ with $\varpi_c\in\mathcal{E}_c$. Define number
		\begin{align}\label{eq43}
			&d:=d\left(\bm{\varpi}\right)=\dfrac{1}{2}\left(\varpi_\infty-\sum_{c\in\Upsilon}\varpi_c\right).
		\end{align}
		We try to find all elements  $\bm{\varpi}$ such that $d$ is a non-negative integer, and retain such elements to perform Step 3.  If there is no such element  $\bm{\varpi}$,  then statement (ii) of Theorem \ref{th6} is impossible.
		\item [\bf{Step 3.}] For each $\bm{\varpi}$ retained from Step 2,  we introduce the rational function
		$$\theta=\dfrac{1}{2}\sum_{c\in\Upsilon}\dfrac{\varpi_c}{z-c}.$$
		Then, we seek a monic polynomial $P$ of degree $d$ defined in \eqref{eq43} such that
		\begin{align*}
			&P'''+3\theta P''+\left(3\theta^2+3\theta'-4r\right)P'+\left(\theta''+3\theta\theta'+\theta^3-4r\theta-2r'\right)P=0.
		\end{align*}
If such polynomial $P$ does not exist for all  elements $\bm{\varpi}$ retained from Step 2,  then statement (ii) of Theorem \ref{th6} is untenable.

Assume that  such a polynomial $P$ exists. Let $\phi=\theta+P'/P$ and $\omega$ be a root of
$$\omega^2-\phi \omega+\left(\dfrac{1}{2}\phi'+\dfrac{1}{2}\phi^2-r\right)=0.$$
Then, $\chi=\exp\left(\int \omega \right)$ is a solution of differential equation $\chi''=r\chi$.
	\end{itemize}
  \end{appendices}


\begin{thebibliography}{222}
\bibitem{MR3743739}
{\sc P.~Acosta-Hum\'{a}nez, M.~Alvarez-Ram\'{\i}rez, and T.~J. Stuchi}, {\em
	Nonintegrability of the {A}rmbruster-{G}uckenheimer-{K}im quartic
	{H}amiltonian through {M}orales-{R}amis theory}, SIAM J. Appl. Dyn. Syst., \textbf{17}
(2018), ~78--96.

\bibitem{MR2784332}
{\sc B.~Andrews and J.~Clutterbuck}, {\em Proof of the fundamental gap
	conjecture}, J. Amer. Math. Soc., \textbf{24} (2011), ~899--916.

\bibitem{MR1218744}
{\sc M.~S. Ashbaugh and R.~D. Benguria}, {\em Eigenvalue ratios for
	{S}turm-{L}iouville operators}, J. Differential Equations, \textbf{103} (1993),~205--219.

\bibitem{MR1081670}
{\sc M.~S. Ashbaugh, E.~M. Harrell, II, and R.~Svirsky}, {\em On minimal and
	maximal eigenvalue gaps and their causes}, Pacific J. Math., \textbf{147} (1991),~1--24.

\bibitem{MR4150221}
{\sc M.~S. Ashbaugh and D.~Kielty}, {\em Spectral gaps of 1-{D} {R}obin
	{S}chr\"{o}dinger operators with single-well potentials}, J. Math. Phys., \textbf{61}
(2020),~091507, 15.

\bibitem{MR0060642}
{\sc P.~F. Byrd and M.~D. Friedman}, {\em Handbook of elliptic integrals for
	engineers and physicists}, Springer-Verlag, 1954.

\bibitem{MR2782621}
{\sc L.~M. Chasman}, {\em An isoperimetric inequality for fundamental tones of
	free plates}, Comm. Math. Phys., \textbf{303} (2011),~421--449.

\bibitem{MR3478937}
{\sc D.~Chen, T.~Zheng, and H.~Yang}, {\em Estimates of the gaps between
	consecutive eigenvalues of {L}aplacian}, Pacific J. Math., \textbf{282} (2016),~293--311.

\bibitem{MR2881964}
{\sc D.-Y. Chen and M.-J. Huang}, {\em Comparison theorems for the eigenvalue
	gap of {S}chr\"{o}dinger operators on the real line}, Ann. Henri
Poincar\'{e}, \textbf{13} (2012),~85--101.

\bibitem{MR4379307}
{\sc H.~Chen, M.~Bhakta, and H.~Hajaiej}, {\em On the bounds of the sum of
	eigenvalues for a {D}irichlet problem involving mixed fractional
	{L}aplacians}, J. Differential Equations, \textbf{317} (2022), ~1--31.

\bibitem{MR3412394}
{\sc H.~Chen and P.~Luo}, {\em Lower bounds of {D}irichlet eigenvalues for some
	degenerate elliptic operators}, Calc. Var. Partial Differential Equations, \textbf{54}
(2015), ~2831--2852.

\bibitem{MR2846268}
{\sc Q.-M. Cheng and G.~Wei}, {\em A lower bound for eigenvalues of a clamped
	plate problem}, Calc. Var. Partial Differential Equations, \textbf{42} (2011),~579--590.

\bibitem{MR639355}
{\sc S.~Y. Cheng and P.~Li}, {\em Heat kernel estimates and lower bound of
	eigenvalues}, Comment. Math. Helv., \textbf{56} (1981),~327--338.

\bibitem{CC77}
{\sc D.~V. Choodnovsky and G.~V. Choodnovsky}, {\em Completely integrable class
	of mechanical systems connected with {K}orteweg-de {V}ries and multicomponent
	{S}chr\"{o}dinger equations}, Lett. Nuovo Cimento (2), \textbf{22} (1978), ~47--51.

\bibitem{MR4075569}
{\sc T.~Combot, A.~J. Maciejewski, and M.~Przybylska}, {\em Bi-homogeneity and
	integrability of rational potentials}, J. Differential Equations, \textbf{268} (2020),~7012--7028.

\bibitem{MR3488539}
{\sc V.~Gol'dshtein and A.~Ukhlov}, {\em On the first eigenvalues of free
	vibrating membranes in conformal regular domains}, Arch. Ration. Mech. Anal.,
\textbf{221} (2016), ~893--915.

\bibitem{MR4447102}
{\sc S.~Guo, G.~Meng, P.~Yan, and M.~Zhang}, {\em Optimal maximal gaps of
	{D}irichlet eigenvalues of {S}turm-{L}iouville operators}, J. Math. Phys., \textbf{63}
(2022), ~Paper No. 072701, 11.

\bibitem{MR4339006}
{\sc J.~Hedhly}, {\em Eigenvalue ratios for vibrating string equations with
	single-well densities}, J. Differential Equations, \textbf{307} (2022),~476--485.

\bibitem{MR1948113}
{\sc M.~Horv\'{a}th}, {\em On the first two eigenvalues of {S}turm-{L}iouville
	operators}, Proc. Amer. Math. Soc., \textbf{131} (2003),~1215--1224.

\bibitem{MR4285920}
{\sc Y.~Ilyasov and N.~Valeev}, {\em Recovery of the nearest potential field
	from the {$m$} observed eigenvalues}, Phys. D, \textbf{426} (2021), ~Paper No.
132985, 6.

\bibitem{MR3912726}
{\sc Y.~S. Ilyasov and N.~F. Valeev}, {\em On nonlinear boundary value problem
	corresponding to {$N$}-dimensional inverse spectral problem}, J. Differential
Equations, \textbf{266} (2019), ~4533--4543.

\bibitem{MR1764944}
{\sc E.~Julliard~Tosel}, {\em Meromorphic parametric non-integrability; the
	inverse square potential}, Arch. Ration. Mech. Anal., \textbf{152} (2000),~187--205.

\bibitem{MR1614731}
{\sc S.~Karaa}, {\em Extremal eigenvalue gaps for the {S}chr\"{o}dinger
	operator with {D}irichlet boundary conditions}, J. Math. Phys., \textbf{39} (1998),~2325--2332.

\bibitem{MR277789}
{\sc T.~Kimura}, {\em On {R}iemann's equations which are solvable by
	quadratures}, Funkcial. Ekvac., \textbf{12} (1969), ~269--281.

\bibitem{MR2299195}
{\sc S.~Kondej and I.~Veseli\'{c}}, {\em Lower bounds on the lowest spectral
	gap of singular potential {H}amiltonians}, Ann. Henri Poincar\'{e}, \textbf{8} (2007),~109--134.

\bibitem{MR839134}
{\sc J.~J. Kovacic}, {\em An algorithm for solving second order linear
	homogeneous differential equations}, J. Symbolic Comput., \textbf{2} (1986),~3--43.

\bibitem{MR1411677}
{\sc V.~V. Kozlov}, {\em Symmetries, topology and resonances in {H}amiltonian
	mechanics}, vol.~31 of Ergebnisse der Mathematik und ihrer Grenzgebiete (3), Springer-Verlag, Berlin,
1996.


\bibitem{MR2682403}
{\sc G.~Kristensson}, {\em Second order differential equations}, Springer, New
York, 2010.


\bibitem{MR1165859}
{\sc P.~Kr\"{o}ger}, {\em Upper bounds for the {N}eumann eigenvalues on a
	bounded domain in {E}uclidean space}, J. Funct. Anal., \textbf{106} (1992),
~353--357.

\bibitem{MR701919}
{\sc P.~Li and S.~T. Yau}, {\em On the {S}chr\"{o}dinger equation and the
	eigenvalue problem}, Comm. Math. Phys., \textbf{88} (1983), ~309--318.

\bibitem{MR573436}
{\sc E.~H. Lieb}, {\em The number of bound states of one-body {S}chroedinger
	operators and the {W}eyl problem}, in Geometry of the {L}aplace operator, Proc. Sympos. Pure Math., XXXVI, Amer. Math. Soc., Providence, R.I.,
1980, ~241--252.

\bibitem{MR2123446}
{\sc A.~J. Maciejewski, M.~Przybylska, and T.~Stachowiak}, {\em
	Non-integrability of {G}ross-{N}eveu systems}, Phys. D, \textbf{201} (2005),~249--267.

\bibitem{MR2831794}
{\sc A.~J. Maciejewski, M.~Przybylska, and A.~V. Tsiganov}, {\em On algebraic
	construction of certain integrable and super-integrable systems}, Phys. D,
\textbf{240} (2011),~1426--1448.

\bibitem{MZ10}
{\sc G.~Meng and M.~R. Zhang}, {\em Continuity in weak topology: first order
	linear systems of {ODE}}, Acta Math. Sin. (Engl. Ser.), \textbf{26} (2010),
~1287--1298.

\bibitem{MR1713573}
{\sc J.~J. Morales~Ruiz}, {\em Differential {G}alois theory and
	non-integrability of {H}amiltonian systems}, vol.~179 of Progress in
Mathematics, Birkh\"{a}user Verlag, Basel, 1999.

\bibitem{MR129219}
{\sc G.~P\'{o}lya}, {\em On the eigenvalues of vibrating membranes}, Proc.
London Math. Soc. (3), \textbf{11 }(1961), ~419--433.

\bibitem{PT87}
{\sc J.~P\"{o}schel and E.~Trubowitz}, {\em Inverse spectral theory}, vol.~130
of Pure and Applied Mathematics, Academic Press, Inc., Boston, MA, 1987.

\bibitem{MR3804721}
{\sc M.~Shibayama}, {\em Non-integrability of the spacial {$n$}-center
	problem}, J. Differential Equations, \textbf{265} (2018), ~2461--2469.

\bibitem{MR829055}
{\sc I.~M. Singer, B.~Wong, S.-T. Yau, and S.~S.-T. Yau}, {\em An estimate of
	the gap of the first two eigenvalues in the {S}chr\"{o}dinger operator}, Ann.
Scuola Norm. Sup. Pisa Cl. Sci. (4), \textbf{12} (1985),~319--333.

\bibitem{MR1383015}
{\sc R.~S. Strichartz}, {\em Estimates for sums of eigenvalues for domains in
	homogeneous spaces}, J. Funct. Anal., \textbf{137} (1996), ~152--190.

\bibitem{MR257592}
{\sc C.~J. Thompson}, {\em On the ratio of consecutive eigenvalues in
	{$N$}-dimensions}, Studies in Appl. Math., \textbf{48} (1969), ~281--283.

\bibitem{TWZ}
{\sc Y.~Tian, Q.~Wei, and M.~Zhang}, {\em On the polynomial integrability of
	the critical systems for optimal eigenvalue gaps}, J. Math. Phys., \textbf{64} (2023),
p.~092701.

\bibitem{MR4197914}
{\sc N.~F. Valeev and Y.~S. Ilyasov}, {\em Inverse spectral problem for
	{S}turm-{L}iouville operator with prescribed partial trace}, Ufa Math. J., \textbf{12}
(2020), ~19--29.

\bibitem{MR1960772}
{\sc M.~van~der Put and M.~F. Singer}, {\em Galois theory of linear
	differential equations}, vol.~328,
Springer-Verlag, Berlin, 2003.

\bibitem{MR2511892}
{\sc H.~Vogt}, {\em A lower bound on the first spectral gap of
	{S}chr\"{o}dinger operators with {K}ato class measures}, Ann. Henri
Poincar\'{e}, \textbf{10} (2009), ~395--414.

\bibitem{WMZ09}
{\sc Q.~Wei, G.~Meng, and M.~Zhang}, {\em Extremal values of eigenvalues of
	{S}turm-{L}iouville operators with potentials in {$L^1$} balls}, J.
Differential Equations, \textbf{247} (2009),~364--400.

\bibitem{MR1511560}
{\sc H.~Weyl}, {\em \"{U}ber gew\"{o}hnliche {D}ifferentialgleichungen mit
	{S}ingularit\"{a}ten und die zugeh\"{o}rigen {E}ntwicklungen
	willk\"{u}rlicher {F}unktionen}, Math. Ann., \textbf{68} (1910),~220--269.

\bibitem{MR0178117}
{\sc E.~T. Whittaker and G.~N. Watson}, {\em A course of modern analysis. {A}n
	introduction to the general theory of infinite processes and of analytic
	functions: with an account of the principal transcendental functions},
Cambridge University Press, New York, 1962.


\bibitem{MR2004534}
{\sc S.~Wiggins}, {\em Introduction to applied nonlinear dynamical systems and
	chaos}, vol.~2 of Texts in Applied Mathematics, Springer-Verlag, New York,
second~ed., 2003.

\bibitem{YZ11}
{\sc P.~Yan and M.~Zhang}, {\em Continuity in weak topology and extremal
	problems of eigenvalues of the {$p$}-{L}aplacian}, Trans. Amer. Math. Soc.,
\textbf{363} (2011),~2003--2028.

\bibitem{YZ12}
{\sc P.~Yan and M.~Zhang}, {\em A survey on extremal
	problems of eigenvalues}, Abstr. Appl. Anal.,  (2012), pp.~Art. ID 670463,
26.

\bibitem{MR923886}
{\sc H.~Yoshida}, {\em A criterion for the nonexistence of an additional
	integral in {H}amiltonian systems with a homogeneous potential}, Phys. D, \textbf{29}
(1987), ~128--142.

\bibitem{MR943701}
{\sc H.~Yoshida}, {\em Nonintegrability of
	the truncated {T}oda lattice {H}amiltonian at any order}, Comm. Math. Phys.,
\textbf{116} (1988), ~529--538.

\bibitem{Zh08}
{\sc M.~Zhang}, {\em Continuity in weak topology: higher order linear systems
	of {ODE}}, Sci. China Ser. A, \textbf{51} (2008),~1036--1058.

\bibitem{Zh09}
{\sc M.~Zhang}, {\em Extremal values of smallest eigenvalues of {H}ill's
	operators with potentials in {$L^1$} balls}, J. Differential Equations, \textbf{246}
(2009), ~4188--4220.
\end{thebibliography}
\end{document}